\documentclass[10pt]{article}
\usepackage{cite,amssymb,graphicx,caption}
\usepackage{amsfonts, amsmath}
\usepackage{enumerate}
\usepackage{epsfig}
\usepackage{xcolor}
\usepackage{psfrag}
\usepackage{hyperref}
\usepackage{babel}
\usepackage{lineno}
\usepackage{rotating}
\usepackage{tcolorbox}
\usepackage[margin=1in]{geometry}
\usepackage{rotating}
\usepackage[colorinlistoftodos]{todonotes}
\usepackage{pgfplots}
\usepackage{stackengine}
\usepackage{authblk}
\usepackage{makecell}
\usepgfplotslibrary{fillbetween}
\pgfplotsset{compat=newest}
\usetikzlibrary{decorations.markings, patterns}

\usepackage{tikz}
\usetikzlibrary{math}
\usepackage{pgfplots}
\usepackage{stackengine}
\usepgfplotslibrary{fillbetween}
\pgfplotsset{compat=newest}
\usetikzlibrary{decorations.markings, patterns}

\newtheorem{lemma}{Lemma}
\newtheorem{prop}{Proposition}
\newtheorem{remark}{Remark}

\definecolor{Bluish}{rgb}{0.,0.,0.5}
\definecolor{Reddish}{rgb}{0.5,0.,0.}
\definecolor{PiVD_color}{RGB}{48,136,163}
\definecolor{ColorDFE}{rgb}{0.035294117647059,0.050980392156863,0.215686274509804}
\definecolor{ColorNME}{rgb}{0.439215686274510,0.737254901960784,0.482352941176471}
\definecolor{ColorNmutE}{rgb}{0,0,0}
\definecolor{ColorCSE}{rgb}{0.188235294117647,0.533333333333333,0.639215686274510}

\def\qed{\hfill $\Box$ \\ \bigskip}
\def\eps{\varepsilon}
\def\rme{\mathrm{e}}
\def\rmi{\mathrm{i}}

\newcommand{\Phit}{\Phi}  
\def\mucrit{\mu_{\textrm{c}}}

\usepackage[colorinlistoftodos]{todonotes}
\usepackage[normalem]{ulem}
\usepackage{color}
\usepackage{comment}
\definecolor{ColorTomasL}{rgb}{0.7,0.2,0.2}

\definecolor{ColorJC}{rgb}{0.1,0.3,0.9}

\definecolor{ColorSanti}{rgb}{0.9,0.,0.4}

\definecolor{ColorJosep}{rgb}{0.1,0.2,1}   

\newcommand{\js}[1]{\textcolor{ColorJosep}{#1}}

\graphicspath{ {./figures/}}

\usepackage{tikz}

\newcommand{\QSmicromu}{\mathrm{QS}_{\mathrm{mic}}^{(\mu)}}
\newcommand{\QSmicromuc}{\mathrm{QS}_{\mathrm{mic}}^{(\mu_c^\ast)}}
\newcommand{\QSmicro}{\mathrm{QS}_{\mathrm{mic}}}

\title{Multiscale feedback drives viral evolution and epidemic dynamics}

\date{\vspace{-5ex}}
\author[1,2,*]{Juan C. Mu\~noz-S\'anchez}
\author[3,4,5,6,*,+]{J. Tom\'as L\'azaro}
\author[5,6]{Josep Sardany\'es}
\author[1,7]{\\Santiago F. Elena}
\affil[1]{Institute for Integrative Systems Biology (I$^2\!$SysBio), CSIC-Universitat de Val\`encia, Paterna, 46980 Val\`encia, Spain}
\affil[2]{Departament de F\'isica Te\`orica, Universitat de Val\`encia, Burjassot, 46100 Val\`encia, Spain}
\affil[3]{Departament de Matem\`atiques, Universitat Polit\`ecnica de Catalunya (UPC), 08028 Barcelona, Spain}
\affil[4]{Institute of Mathematics, UPC-BarcelonaTech (IMTech), 08028 Barcelona, Spain}
\affil[5]{Centre de Recerca Matem\`atica (CRM), Cerdanyola del Vall\`es, 08193 Barcelona, Spain}
\affil[6]{Dynamical Systems and Computational Virology, CSIC Associated Unit CRM-I$^2\!$SysBio, Spain}
\affil[7]{Santa Fe Institute, Santa Fe, NM 87501, USA}
\affil[*]{Equal contribution}
\affil[+]{Correspondence: jose.tomas.lazaro@upc.edu}

\date{\today}

\begin{document}

\maketitle


\section*{Abstract}

We introduce a minimal multiscale framework that links within‑host virus dynamics to population‑level SIRS epidemiology through explicit, bidirectional coupling. At the microscopic layer, a two variant quasispecies (master and mutant genomes with packaged virions) evolves on a fast timescale. At the macroscopic layer, two infectious classes (master- and mutant-infected), susceptible, recovered, and deceased individuals evolve slowly. The two scales are connected through transmission rates that depend on instantaneous virion abundance and through prevalence‑weighted effective replication rates. Exploiting the timescale separation, we formalize a coarse‑grained slow-fast closure: the genome-virion subsystem rapidly relaxes to quasi-steady states that parameterize time-varying transmission in the slow epidemiological system. This yields an integrated expression for the basic reproduction number and sharp inequalities that delineate coexistence \emph{versus} exclusion. A key prediction is a context-dependent error threshold that shifts with the prevalence ratio, enabling transient pseudo-error catastrophes driven by epidemic composition rather than intrinsic fidelity. Linearization reveals parameter regions with damped oscillations arising solely from the microscopic-macroscopic feedback. Two illustrative extremes bracket the model’s behavior: an avirulent strongly immunizing strain that benignly replaces the master, and a hypervirulent weakly immunizing that self-limits via host depletion and collapses transmission. This framework yields testable signatures linking viral load, incidence, and within‑host composition.

\section{Introduction}
Viruses infect every form of life, from viruses and microbes to plants and animals. Most coexist with their hosts, but spillovers into new hosts can cause disease and reduce host fitness. At the same time, viruses can drive host innovations (\textit{e.g.},~\cite{Mi2000,ShapiroTurner2018}). Viral biology is studied across multiple scales: molecular virologists dissect replication, pathogenesis, and antiviral defenses within cells; clinicians focus on individuals; and ecologists and epidemiologists track incidence, transmission, and ecosystem-level impacts. Accordingly, what constitutes the ``host" (cell or tissue, individual or population and ecosystem) and the operative ``viral unit" (mutant swarms within individuals, infected hosts in ecological studies, or evolutionary lineages in phylogeography) shifts with scale. Yet the ultimate host is the cell, where viral genes are expressed, replication organelles form, genomes and proteins are synthesized, host defenses are countered, particles assemble, and infection spreads to neighboring cells, tissues, new hosts, and ultimately through populations and ecosystems. Selective pressures also differ across scales: rapid replication may confer an advantage within hosts but prove suboptimal for between‑host transmission \cite{Doumayrou2012}. Thus, although the qualitative link between within-host infection dynamics and population‑level transmission is widely acknowledged, a comprehensive, quantitative framework that integrates these scales remains lacking.

As complex dynamical systems, viruses can exhibit phase transitions involving abrupt changes in dynamical behavior or internal structure~\cite{Sole2021}.  At the molecular level, critical phenomena are exemplified by the error threshold associated to highly mutagenic replication of viral genomes~\cite{Eigen1971,Eigen1989}. At the epidemiological level, the spread of viral diseases across heterogeneous contact networks displays complex dynamical behavior and phase transitions associated with the existence of highly connected hubs \cite{Barthelemy2005,BalcanVespignani2011}. So far, models connecting phase transitions across biological scales deserve further exploration. In particular, no comprehensive framework has yet been developed to investigate whether transitions at one level may result from the dynamical properties of the lower levels.

There have been few attempts to model multiscale selection in viruses. Existing studies have treated between-host transmissions as a function of within-host replication parameters \cite{Coombs2007,Mideo2008,Feng2012,Scholle2013,ShinMacCarthy2016,Doekes2017,Dorratoltaj2017}, or tissue and organ colonization as an extension of within-cell replication processes and their interaction with host factors \cite{SardanyesElena2011,Heldt2013,Kumberger2016}.  Unfortunately, the former models usually neglect the inherent within-host complexity, while the latter rarely extend beyond individual tissues or single hosts. Modeling multiscale processes in full mechanistic detail is therefore unfeasible. One alternative is to extract the essential features of lower-scale models to embed them into higher-scale descriptions~\cite{Mideo2008}. An approach that has been successfully taken is to separate timescales, which enables the construction of effective reduced models that operate at distinct biological scales.

Predicting when a viral mutant genotype will rise to dominance and potentially trigger an epidemic is a critical challenge. This problem inherently involves two scales: a microscopic scale, at which a mutant genome is generated from a master sequence, and a macroscopic, or population scale, at which the new virus must sustain a transmission chain to remain in circulation and outcompete existing strains. At the microscopic level, quasispecies theory provides a dynamical description of replicator evolution under mutation--selection~\cite{Eigen1971,Eigen1977}. Within this framework, a quasispecies is understood as a structured distribution of closely related variants maintained by the balance between mutation and selection, rather than a single dominant genotype~\cite{Nowak1992}. In this sense, selection acts on the population as a whole, shaping the collective composition of variants rather than favoring an individual sequence in isolation. A fundamental result of quasispecies theory is the existence of an error threshold, beyond which mutation dominates over selection and the population delocalizes in sequence space, preventing the stable maintenance of genetic information~\cite{Eigen1971,Eigen1977,Eigen1989}. This threshold has been experimentally identified in RNA viruses~\cite{Crotty2001} and has been proposed to operate in hepatitis C virus-infected patients through replicative, rather than mutational, thresholds~\cite{Sole2006}. At the population level, a range of epidemiological models describes viral transmission in terms of biologically meaningful parameters, with the Susceptible-Infectious-Recovered (SIR) model and its variants being among the most fundamental. 

These two scales are deeply interconnected. Evolutionary pressures acting on viral populations may depend on macroscopic factors such as host availability, population immunity, transmission opportunities, or behavioral changes, while the epidemiological success of a variant depends on microscopic properties such as replication rate, mutation rate, and within-host viral production. The COVID-19 pandemic provides a clear example of this coupling, with co-circulating SARS-CoV-2 variants emerging through random mutation and spreading under the combined action of biological, epidemiological, and social selection pressures~\cite{Carabelli2023}. Understanding variant emergence and competition therefore requires a framework able to connect intracellular viral dynamics with population-level transmission.

In this work, we propose a multiscale modelling approach that integrates known intracellular replication dynamics with macroscopic population-level transmission. Our objective is to better understand the characteristics a mutant virus must possess to become dominant, potentially replacing the master sequence entirely. The link between both scales has been explicitly made by expressing the macroscopic transmission rates as a function of the number of viral particles produced in infected cells. In addition, we explore how population level processes can drive the extinction of specific viral strains and how intracellular dynamics can trigger nontrivial epidemic events.
 
Among our findings, we present a formulation for the basic reproduction number, $\mathcal{R}_0$, that incorporates parameters from both scales, thereby providing a more comprehensive measure of viral fitness. Among the many scenarios that can be analyzed within our framework, we illustrate two extremes along a continuum: an avirulent virus with rapid recovery that induces temporal immunity (\emph{i.e.}, a vaccine-like virus) and a hypervirulent virus with slow recovery (\emph{i.e.}, a burnout virus). In the first scenario, the master sequence is replaced by the avirulent mutant, which persists at a high incidence in the host population. In the second scenario, all infected individuals eventually die, extinguishing the infection, leaving only susceptible individuals in the population. Together, these examples show how coupling microscopic mutation--selection dynamics with macroscopic transmission can generate qualitatively different epidemic regimes and provide a minimal framework for studying viral emergence across biological scales.

\section{Multiscale mathematical model}
 
\begin{figure}[!ht]
    \vspace{-0.2cm}
    \centering
    \begin{minipage}{0.38\textwidth}
        \raggedright
        \includegraphics[scale=0.3]{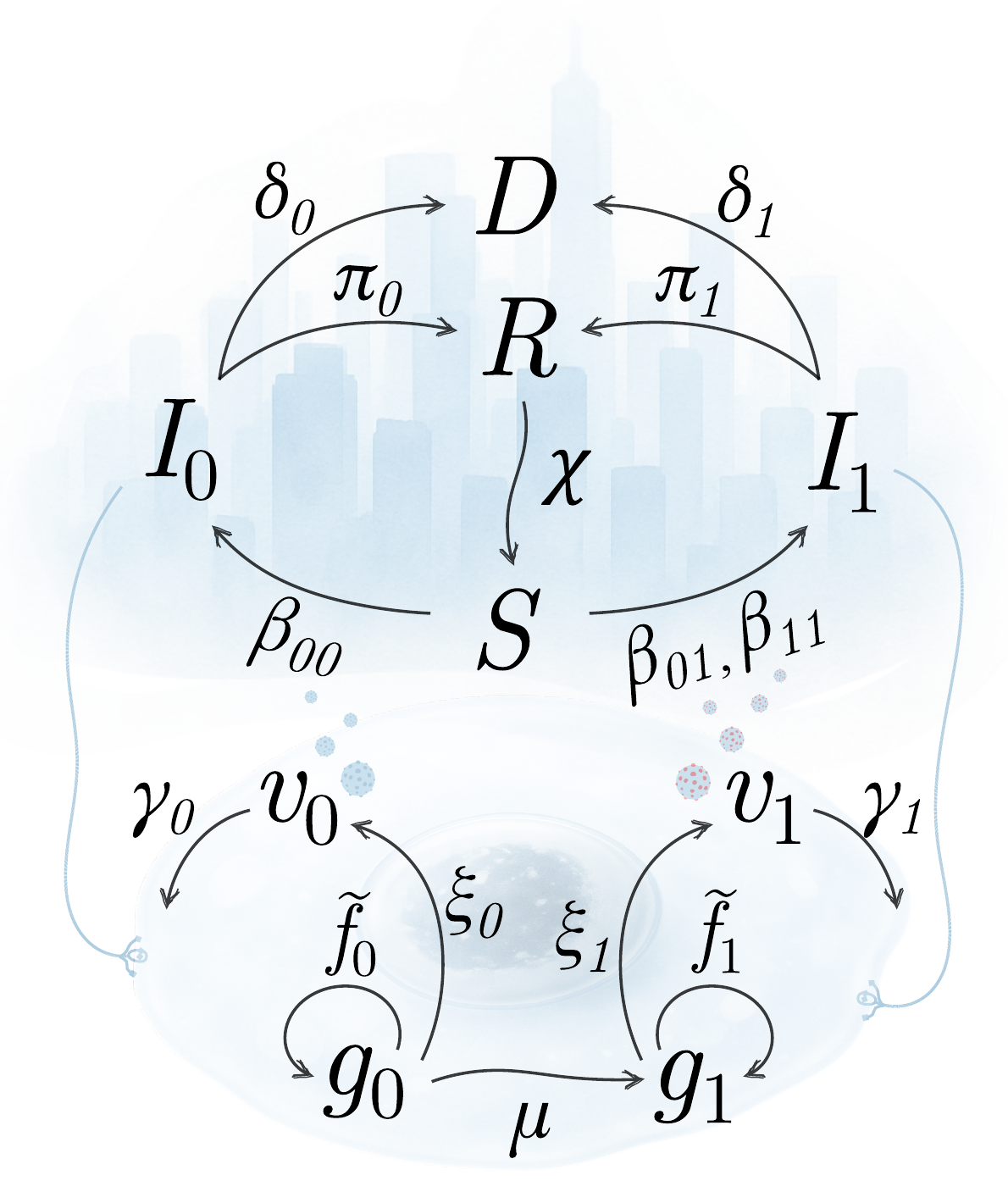}
    \end{minipage}
    \hfill
    \begin{minipage}{0.61\textwidth}
        \raggedright
        \vspace{+0.75cm}
        \begin{tabular}{c l l}
            \hline
            Parameter    & Description & Stoichiometry\\
            \hline\noalign{\vskip 4pt}
            $\delta_{j}$   & Mortality rate & $I_j\xrightarrow{\delta_j} D \quad j=0,1$\\
            $\pi_{j}$      & Recovery rate & $I_i\xrightarrow{\pi_j} R \quad j=0,1$\\
            $\chi$           & Waning immunity rate & $R \xrightarrow{\chi} S$\\
            $\beta_{ij}$     & Transmission rates & \makecell[tl]{$S+I_i \xrightarrow{\beta_{ij}} I_i+I_j,$\\  \hspace{34pt} {\small $ij\in\{00,01,11\}$}}\\
            $\tilde{f}_{j}$        & Effective fitness rate  & $g_j\xrightarrow{\tilde{f}_j}g_j+g_j\quad j=0,1$\\
            $\mu$            & Mutation probability & $\mu g_0 \xrightarrow{\tilde{f}_{0}} g_1$\\
            $\xi_{j}$      & Virions production rate & $g_j \xrightarrow{\xi_{j}} v_j\quad j=0,1$\\
            $\gamma_{j}$   & Virions decay rate& $v_j \xrightarrow{\gamma_f}\varnothing\quad j=0,1$\\[3pt]
            \hline
        \end{tabular}
    \end{minipage}
    \caption{Schematic representation of the multiscale model (left) and summary of the associated parameters (right). Microscopic variables include the master and mutated genomes ($g_0$, $g_1$) and their corresponding virions ($v_0$, $v_1$), while macroscopic variables comprise susceptible ($S$), infected ($I_0$, $I_1$), recovered ($R$), and deceased ($D$) individuals.}
    \label{fig:schematicModel}
\end{figure}

We propose a minimal multiscale model that links within‑host and population dynamics: a quasispecies module describing mutation from a master to a single dominant mutant genome within an average host is coupled to a two‑strain SIRS system that distinguishes individuals infected by each variant. Here, \emph{multiscale} denotes the coupling of processes across distinct levels of biological organization (within-host and between-host), which is translated in our model into a separation of characteristic timescales. This usage is well established in infectious disease modeling and evolutionary epidemiology~\cite{Mideo2008,Alizon2011,Handel2015}. For analytical clarity, mutation is assumed to be unidirectional (no reverse mutation or lineage diversification). This is a common assumption in quasispecies models: the probability of backward mutations is extremely low due to the enormous size of the sequence space~\cite{Sole2003}. The model comprises two layers and two time scales (Figure~\ref{fig:schematicModel}): a fast microscopic layer tracking genome and virion abundances for the master $(g_0,v_0)$ and for the averaged mutant $(g_1,v_1)$, and a slow macroscopic SIRS layer with infectious classes $I_0$ and $I_1$; no inflow of susceptible individuals is included. The ratio between their respective characteristic timescales is set by the parameter $\varepsilon$. This is a classical multiscale approach and has been extensively employed in biological modeling to describe systems with fast and slow dynamics~\cite{KeenerSneyd1998}. The layers are bidirectionally coupled: transmission rates in the SIRS module depend on the virion distribution, while key microscopic parameters are modulated by the current epidemiological state. Despite its idealization, this framework captures the emergence, competition, and replacement of dominant variants and clarifies conditions that favor the rise of new strains. The detailed models at each level are:

\begin{itemize}
\item\textbf{Microscopic level.}
At the genome's level the model follows a variant of the Swetina-Schuster quasispecies approach~\cite{Swetina1982}:
\begin{eqnarray}
\eps \dot{g}_0 &=& f_0 (1-\mu) \nu_0 g_0 - \Phit g_0, \label{eq:g:0}\\
\eps \dot{g}_1 &=& f_0 \mu \nu_0 g_0 + f_1 \nu_1 g_1 - \Phit g_1, \label{eq:g:1}
\end{eqnarray}
where $g_0$ and $g_1$ are the master and the pool of mutant genomes, respectively, $f_0>f_1>0$ are their fitness rates, and  $0\leq \mu \leq 1$ is the mutation probability. The parameter $\eps>0$ is assumed to be small\footnote{Intracellular viral processes typically occur within a few hours~\cite{Baccam2006,Kakizoe2015}, whereas epidemiological processes such as transmission and recovery evolve over several days~\cite{Chan2025,CDCInfluenza2025}; a separation of roughly one to two orders of magnitude, consistent with the representative value used in the simulations.}. Moreover, being $I_0$ and $I_1$ the prevalence of two different types of infectious individuals, we denote by
\begin{equation}
\nu_0 = \nu_0(I_0,I_1)=\frac{I_0}{I_0+I_1}\quad \textrm{and} \quad \nu_1 = \nu_1(I_0,I_1)=\frac{I_1}{I_0+I_1}  \qquad \textrm{if} \quad (I_0,I_1) \neq (0,0),
\label{def:nu}
\end{equation}
and $0$ otherwise, the relative prevalence of $I_i$ with respect to the total infected individuals $I_0+I_1$. 
Since, in our case, all the infectious individuals arise from these two types, it follows that $\nu_1 = 1 - \nu_0$, except when $I_0=I_1=0$, in which $\nu_0=\nu_1=0$.

The function $\Phit$ in~\eqref{eq:g:0}-\eqref{eq:g:1} is the standard outflow term introduced to keep the total population of genomes, $g_0+g_1$, constant. In this scenario this is given by
\[
\Phi=f_0 \nu_0 \frac{g_0}{g_0+g_1} + f_1 \nu_1 \frac{g_1}{g_0+g_1}.
\]
If we take $g_0,g_1$ such that $g_0(0)+g_1(0)=1$ (so, we will talk about frequencies) this function becomes
\[
\Phi=f_0 \nu_0 g_0 + f_1 \nu_1 g_1
\]
and hence $g_0(t)+g_1(t)=1$ $\forall t\geq 0$. This will be assumed along this work. Whenever both infected populations $I_0, I_1$ simultaneously vanish (total virus extinction), the system~\eqref{eq:g:0}-\eqref{eq:g:1} becomes
$\eps \dot{g}_0 = 0$, $\eps \dot{g}_1 =0$.

\bigskip

It is well-known in quasispecies theory, the existence of the so-called
critical mutation driving the system into error catastrophe, a value for $\mu$ beyond which the master sequence disappears and only mutants persist~\cite{Eigen1971,Eigen1989,Sole2021}. In our case, it is given~\footnote{It suffices to define the effective fitnesses $\tilde{f}_0= f_0 \nu_0$ and $\tilde{f}_1=f_1 \nu_1 = f_1 (1-\nu_0)$ and to use the classical formula $\mucrit = 1 - \tilde{f}_1/\tilde{f}_0$.}  by
\begin{equation}
\mucrit = 1 - \frac{f_1}{f_0} \left( \frac{1}{\nu_0} - 1 \right),
\label{def:genom:critical:mutation:rate}    
\end{equation}
only defined for $I_0>0$ (otherwise $g_0$ cannot replicate).
The terms $\nu_0$ and $\nu_1$ connect the microscopic to the macroscopic layers by incorporating the assumption that infected individuals become a necessary substrate for the replication of viral genomes. The fitness rates are modulated by the fraction of infected individuals of each type.

Concerning the virion's layer, we consider
\begin{eqnarray}
\eps\dot{v}_0 &=& \xi_0 g_0 - \gamma_0 v_0, \label{eq:v:0} \\
\eps\dot{v}_1 &=& \xi_1 g_1 - \gamma_1 v_1,  \label{eq:v:1}  
\end{eqnarray}
where $v_0$ and $v_1$ are the corresponding virions for the genomes $g_0$ and $g_1$, respectively,
$\xi_{0,1}>0$ denote their encapsidation constants, and $\gamma_{0,1}>0$ their degradation rates.

Notice, for the genome's and virion's systems, the inclusion  of the small parameter $\eps>0$ which makes the microscopic level to evolve much faster in time than its counterpart macroscopic one (see below). This leads to a two-timescales model.

\item \textbf{Macroscopic level.} This is defined as a SIRS epidemiological model with feedback  with the microscopic layer.
The variables $S$ stands for susceptible individuals, $I_j$ for infected\footnote{The subscript of an infected compartment denotes the strain that initiated the infection in that individual, and, for $I_0$, not necessarily the strain that may be transmitted during subsequent contacts.} by virions' type $v_j$ ($j=0,1$), $R$ for recovered and $D$ for dead. The model does not include co-infection by the master and mutant strains, nor any conversion of infected individuals between the classes $I_0$ and $I_1$.
Precisely, this model reads
\begin{eqnarray}
\dot{S} &=& - \left( \beta_{00} I_0  + \beta_{01}  I_0 + \beta_{11} I_1 \right) S +   \chi R, \label{eq:sirs:S} \\
\dot{I}_0 &=& \beta_{00} I_0 S - \left( \pi_0 + \delta_0 \right) I_0,
\label{eq:sirs:I0} \\
\dot{I}_1 &=& \left( \beta_{01} I_0 + \beta_{11} I_1 \right) S - \left( \pi_1 + \delta_1 \right) I_1, \label{eq:sirs:I1} \\
\dot{R} &=& \pi_0 I_0 + \pi_1 I_1 - \chi R,  \label{eq:sirs:R} \\
\dot{D} &=& \delta_0 I_0 + \delta_1 I_1  \label{eq:sirs:D} 
\end{eqnarray}
where $\delta_{0,1}$ denote the viral-induced mortality rates, $\pi_{0,1}$ the recovery rates, and the infection rates $\beta_{00}, \beta_{01}, \beta_{11}$ that will connect the micro- and macroscopic levels (see below) being functions of the viral composition, which itself depends on the epidemiological state. The parameter $\chi>0$ is the waning immunity rate; it is assumed (for simplicity) to be the same for both type of viruses (master and its mutant). The no-presence of the small parameter $\eps$ makes this macroscopic system slow in comparison with the microscopic level.

Observe that $\dot{S}+\dot{I}_0+\dot{I}_1+\dot{R}+\dot{D}=0$ and so $S(t)+I_0(t)+I_1(t)+R(t)+D(t)$ is a first integral of the system~\eqref{eq:sirs:S}-\eqref{eq:sirs:D}. This implies that $S(t)+I_0(t)+I_1(t)+R(t)+D(t)=S(0)+I_0(0)+I_1(0)+R(0)+D(0)$ for any $t\geq 0$. It is not a loss of generality to assume this value equal to $1$ and so consider the variables $S,I_0,I_1,R$, and $D$ to represent fractions of the total (invariant) population.
\end{itemize}

\renewcommand{\arraystretch}{1.1}

\subsection{Multiscale coupling between microscopic and macroscopic dynamics}

The coupling between microscopic and macroscopic levels in our model arises through the construction of key parameters that bridge the two scales. We propose the following bidirectional connections:
\begin{enumerate}
\item{\textbf{From microscopic to macroscopic dynamics via transmission rates.}} A connection from the intra-host to the population level is established through the transmission rates, denoted by $\beta_{ij}(v_j(t))$. 
These rates quantify the efficiency with which an individual infected with variant $i$ transmits the infection to a susceptible individual, who subsequently develops an infection caused by variant $j$. Crucially, we assume that these rates depend on the instantaneous concentration of virions, $v_j(t)$, associated with each viral variant~\footnote{Recall than within $I_0$ individuals, $g_0$ replicates generating $g_1$. Therefore, both virion types can be produced. In contrast, within $I_1$ individuals only $g_1$ replicates, given our assumption of no back mutation.}. This dependence is modeled using a hyperbolic-like function, which captures the saturating nature of transmission dynamics with increasing viral load.

To account for the different infection pathways, we define three \footnote{Notice that 
$\beta_{10} = 0$, since mutant-infected individuals ($I_1$) do not carry master virions, $v_0^{(1)} = 0$, and thus cannot generate new infections of the master type.} distinct transmission rates:
\begin{itemize}
    \item $\beta_{00}(v_0)$: transmission from a master-infected individual leading to a new master infection.
    \item $\beta_{01}(v_1)$: transmission from a master-infected individual leading to a new mutant infection.
    \item $\beta_{11}(v_1)$: transmission from a mutant-infected individual leading to a new mutant infection.
\end{itemize}
To estimate the effective availability of mutant virions $v_1$ within the population, we assume a proportional distribution of $v_1$ between the two types of infected individuals. Specifically, we define:
\begin{equation*}
    v_1^{(j)} = v_1 \cdot \frac{I_j}{I_0 + I_1} = v_1 \nu_j
\end{equation*}
This assumption reflects the idea that the total mutant viral load is partitioned among infected individuals proportionally to their prevalence. Furthermore, these expressions then inform the variant-specific transmission terms used in the population-level model.

The transmission rates are assumed to follow a saturation function of the form:
\[
\beta_{00}= \beta_{00}(v_0)=\frac{a_0 v_0}{b_0 + v_0}, \quad
\beta_{01}=\beta_{01}(v_1)=\frac{a_1 \nu_0 v_1}{b_1 + \nu_0 v_1 }, \quad
\beta_{11}=\beta_{11}(v_1)=\frac{a_1 \nu_1 v_1}{b_1 + \nu_1 v_1} =
\frac{a_1 (1-\nu_0) v_1}{b_1 + (1-\nu_0)v_1},
\]
with constants $a_{0,1}\geq 0$, $b_{0,1}>0$, and $\nu_0=\nu(I_0,I_1)$ as defined in~\eqref{def:nu}. Biologically, $a_{0,1}$ represents the contagion potential (\textit{i.e.}, the maximum value of the transmission rate function), and $b_{0,1}$ represents the characteristic viral load at which the contagion rate reaches half of its maximal value. The terms $\nu_0v_1$ and $\nu_1v_1$ in $\beta_{01}, \beta_{11}$ represent the fractions of mutant virions produced within master-infected individuals and mutant-infected individuals, respectively. This splitting is a mean-field closure: instead of introducing separate microscopic viral populations for each infected class, their effective contributions are represented by averaged quantities weighted by the prevalence ratios $\nu_0$ and $\nu_1$. Accordingly, the dependence of the transmission coefficients on the epidemiological state should be interpreted as an aggregate approximation of viral composition.

\item{\textbf{Macroscopic to microscopic levels via prevalence.}}
Conversely, to represent how population-level dynamics influence viral evolution at the intracellular scale, we modulate the replicative fitness of each variant based on the relative abundance of host types in the population. Since infected individuals act as the primary replication environment for the virus, the availability of hosts directly affects the success of the variants. Although mutant virions may replicate to some extent in individuals originally infected with the master variant, for simplicity, we assume that each viral variant replicates predominantly within its corresponding host type. Thus, the replicative fitness parameters are defined as functions of the current population state:
\begin{equation*}
    \tilde{f}_0(t) = f_0 \cdot \frac{I_0(t)}{I_0(t) + I_1(t)} = \nu_0(t) f_0,
    \qquad
    \tilde{f}_1(t) = f_1 \cdot \frac{I_1(t)}{I_0(t) + I_1(t)} = \nu_1(t) f_1 = (1-\nu_0(t)) f_1,
\end{equation*}
where $f_0$ and $f_1$ are the baseline fitness values for the master and mutant variants, respectively. These bidirectional links enable the model to capture feedback loops between within-host dynamics and between-host transmission, offering insights into how selective pressures at one scale shape dynamics at the other. 

It is important to note that, in this framework, the term ``multiscale" refers not only to interactions between biological levels (intracellular and population) but also to the separation of time scales between them. Specifically, we assume that microscopic (within-host) dynamics occur on a much faster timescale than  macroscopic (population-level) dynamics.
\end{enumerate}
Regarding notation, from this point onward we denote $\nu = \nu_0$ and, consequently, $\nu_1 = 1 - \nu_0 = 1 - \nu$. This convention fails only in the special case $I_0 = I_1 = 0$, for which $\nu_0 = \nu_1 = 0$. This case is explicitly addressed at the appropriate point in the analysis.

\section{Results and discussion}

\subsection{Equilibrium points}
\label{sec:eqPoints}

We begin by computing the equilibrium points of systems~\eqref{eq:g:0}-\eqref{eq:g:1},~\eqref{eq:v:0}-\eqref{eq:v:1}~and~\eqref{eq:sirs:S}-\eqref{eq:sirs:D}. Together with other invariant objects, these equilibria constitute the dynamical skeleton of the system and fully determine its dynamics. Their analysis is therefore essential for understanding the possible long-term regimes of the model, identifying which viral strains may persist or disappear, and clarifying how the coupling between microscopic mutation dynamics and macroscopic epidemiological processes shapes the global behaviour of the system. Henceforth, we assume that the virions encapsidation and degradation rates are positive, \emph{i.e.} $\xi_{0,1}>0$ and $\gamma_{0,1}>0$,
and that $\mu \in (0,1]$ (\textit{i.e.}, there is always $g_1$ production due to errors in $g_0$ replication).

The different classes of equilibrium points are classified according to the type (if any) of infected individuals that remain:
\begin{itemize}
\item \textbf{DFE}: disease-free equilibrium, when $I_0=I_1=0$.
\item \textbf{NME}: equilibrium with no master-infected individual\js{s}, $I_0=0$, $I_1> 0$.
\item \textbf{NmutE}: equilibrium with no mutant-infected individuals, $I_0> 0$ and $I_1=0$;
\item \textbf{CSE}: co-circulating strains equilibrium, with $I_0> 0$, $I_1> 0$.
\end{itemize}
Notice this is a macroscopic-level criterion, defined from an epidemiological perspective. The explicit expressions of such equilibria are provided in the following propositions. To ease the reading, their proofs have been deferred to Appendix~\ref{sec:proofs:eqPoints}.

We also impose the following assumptions on model parameters to restrict the range of feasible scenarios:
\begin{itemize}
\item \textbf{Microscopic system:} $\mu \in (0,1]$ and $f_0 > f_1 > 0$ (genomes); $\xi_{0,1} \geq 0$ and $\gamma_{0,1} > 0$ (virions).
\item \textbf{Macroscopic system:} $a_{0,1} \geq 0$, $b_{0,1} > 0$, $\delta_{0,1} \geq 0$, $\pi_{0,1} \geq 0$, and $\chi \geq 0$ (individuals).
\end{itemize}

The link between the relative prevalence $\nu(I_0,I_1)$ and the critical mutation probability $\mu_c$ derives into a close connection 
between the microscopic and macroscopic equilibrium points. 
Namely, let us assume $(I_0^*,I_1^*)\ne (0,0)$ being 
$I$-equilibrium infected populations of system~\eqref{eq:sirs:S}-\eqref{eq:sirs:D}. Denote by $\nu^*=\nu_0^*= \nu(I_0^*,I_1^*)\in (0,1)$, as defined in~\eqref{def:nu}, the corresponding relative prevalence of $I_0^*$ with respect to $I_0^*+I_1^*$, and by $\mucrit^*$ its critical mutation rate (see~\eqref{def:genom:critical:mutation:rate}):
\begin{equation}
\mucrit^* = 1 - \frac{f_1}{f_0} \left( \frac{1}{\nu^*} -1 \right).
\label{def:mucrit_star}
\end{equation}
Recall that this value $\mucrit^*$ determines the survival of the master strain or the dominance by the mutant. Their steady states are strongly related to the difference between the current mutation rate $\mu$ and this critical threshold $\mucrit^*$. 
The following lemma specifies the form of these microscopic equilibrium points. It will be taken into account to determine the complete sets of equilibria discussed afterwards.

\begin{lemma}[macroscopic to microscopic equilibrium points]
\label{lemma:micro}
Let $I_0^*$, $I_1^*$ be equilibrium points of the macroscopic system above. Then, their associated genome-virion equilibrium points take the following form:
\[
\begin{array}{lll}
\QSmicromuc: & \ (g_0^*, g_1^* \, ; v_0^*, v_1^*)  =
\left(  1 - \dfrac{\mu}{\mucrit^*} , \dfrac{\mu}{\mucrit^*} \, ;
\dfrac{\xi_0}{\gamma_0}\, g_0^*, \dfrac{\xi_1}{\gamma_1}\, g_1^*
\right) &\qquad \textrm{if $0 < \mu<\mucrit^*$}  \\[2.5ex]
 & \ (g_0^*, g_1^* \, ; v_0^*, v_1^*)  =
\left( 0, 1 \, ; 0, \dfrac{\xi_1}{\gamma_1} \right) &\qquad \textrm{if $\mucrit^* \leq \mu \leq 1$},
\end{array}
\]
\end{lemma}


\vskip1cm

The full set of equilibrium points of the coupled system (that is, the points of type \textbf{DFE}, \textbf{NME},\textbf{NmutE}, and \textbf{CSE} defined above) is characterized in the following propositions. Their 
microscopic and macroscopic expressions are determined and necessary feasibility conditions are explicitly stated.

\bigskip

The first result refers to the existence and form of the \textbf{DFE}-points. These equilibrium points correspond to the lack of infection in the macroscopic system: no individuals remain infected by either the master or the mutant strain. Recall that in this case, $\nu_0^*=\nu_1^*=0$.

\begin{prop}[disease-free equilibrium (\textbf{DFE}) points]
\label{prop:DFE}
Let us assume an equilibrium point satisfying that $I_0^*=I_1^*=0$. Then, the microscopic equilibria take the form:
\[
\QSmicro: (g_0,g_1,v_0,v_1) = \left( g_0, 1-g_0, \frac{\xi_0}{\gamma_0} g_0, \frac{\xi_1}{\gamma_1} (1-g_0) \right), 
\]
for arbitrary $g_0,g_1$ such that $g_0+g_1=1$. Concerning the macroscopic equilibrium point, it is of the form
\begin{equation*}
(S, I_0^*, I_1^*, R, D) \times (g_0, g_1 \, ; v_0, v_1)=
(S, 0,0,R, D) \times \QSmicro ,
\end{equation*}
where $R=0$ if there is no permanent immunity (i.e. $\chi>0$) or $R>0$, otherwise. In all cases, $S+R+D=1.$
\end{prop}

\bigskip

As a general notation rule along this section, concrete equilibrium values will be accompanied by an asterisk, for instance $I_j^*$, $g_j^*, \ldots$; while variables without it, $R, D, S, \ldots$, will correspond to arbitrary values, restricted always, to the total ``mass" being equal to $1$.

The second type of equilibrium points, denoted by \textbf{NME}, describes a situation in which the mutant population dominates at the macroscopic level, with no presence of the master-infected population. As shown in the following proposition, such equilibria are feasible only if mutant-induced mortality vanishes, so that individuals in $I_1$ do not die. The resulting macroscopic and microscopic steady-state values are then determined by the waning immunity rate $\chi$ and the transmission rate $\beta_{11}$.

\begin{prop}[no master equilibrium (NME) points]
\label{prop:NME}
Equilibrium points with no master-infected individuals and surviving mutant-infected individuals, \emph{i.e.} $I_0^*=0, I_1^*>0$, exist if and only if there is no mutant induced mortality, $\delta_1=0$.
Moreover, the microscopic system at equilibrium falls into one of the following two cases:
\[
\QSmicro^{(0)}: (g_0^*,g_1^*,v_0^*,v_1^*)=\left(1 , 0, \dfrac{\xi_0}{\gamma_0},0\right)
\qquad
\textrm{or}
\qquad
\QSmicro^{(1)}: (g_0^*,g_1^*,v_0^*,v_1^*)=\left( 0, 1 , 0, \dfrac{\xi_1}{\gamma_1}\right).
\]
The complete equilibrium set must belong to one of the following cases:
\begin{itemize}
\item[(i)] Case $\chi>0$, \textbf{waning immunity}: they are two possible scenarios, according to the value of the mutant transmission rate $\beta_{11}$.
\begin{itemize}
\item[(i$_1$)] If $\beta_{11}>0$, they are given by
\[
(S^*, I_0^*, I_1^*, R^*, D^*) \times (g_0^*, g_1^* \, ; v_0^*, v_1^*)=\left(\frac{\pi_1}{\beta_{11}},0,I_1^*,\frac{\pi_1}{\chi}I_1^*,D^* \right)\times  \QSmicro^{(1)} ,
\]
\item[(i$_2$)] If $\beta_{11}=0$, then necessarily the mutant recovery rate must also vanish, $\pi_1=0$, and the equilibria are of the form
\[
(S, I_0^*, I_1^*, R^*, D) \times (g_0^*, g_1^* \, ; v_0^*, v_1^*)= \left(S,0,I_1^*,0,D \right)\times
\left\{ \QSmicro^{(0)}\ , \ \QSmicro^{(1)} \ \textrm{iff $a_1=0$ or $\xi_1=0$} \right\},
\]
with arbitrary $S$ and $D$.
\end{itemize}
\item[(ii)] Case $\chi=0$, \textbf{permanent immunity}: in this situation, equilibria exist only if there is no mutant recovery rate, $\pi_1=0$. Like in the previous case $(i)$, their form depend on the transmission rate $\beta_{11}$. Namely,
    \begin{itemize}
    \item[(ii$_1$)] If $\beta_{11}>0$: 
    \[
    (S^*, I_0^*, I_1^*, R, D) \times (g_0^*, g_1^* \, ; v_0^*, v_1^*)=\left(0,0,I_1^*,R, D \right)\times \left\{ \QSmicro^{(0)} \ , \  \QSmicro^{(1)} \right\},
    \]
    with arbitrary $R$ and $D$. 
    \item[(ii$_2$)] If $\beta_{11}=0$: 
    \[
    (S, I_0^*, I_1^*, R, D) \times (g_0^*, g_1^* \, ; v_0^*, v_1^*)=\left(S,0,I_1,R,D \right)\times
    \left\{ \QSmicro^{(0)} \ , \ \QSmicro^{(1)} \ \textrm{iff $a_1=0$ or $\xi_1=0$} \right\},
    \]
    with arbitrary $S$, $R$, and $D$. 
\end{itemize}
\end{itemize}
Recall that, in all cases, the sum of the variables must be equal to $1$.
\end{prop}

\bigskip

The case of equilibrium points of type \textbf{NmutE} shows, in some sense, the symmetric situation to the latter. In this framework, the master strain infection dominates (at macroscopic level) and mutant infected individuals disappear. As in the previous proposition, a key point is the necessary and sufficient condition for the master strain of not being lethal. The variety of possible macroscopic steady states is given by the waning immunity rate $\chi$ and the transmission rates $\beta_{00},\beta_{01}$. On the other side, the microscopic equilibrium depends on the value of the mutation probability $\mu\in (0,1]$. 

\begin{prop}[no mutant equilibrium (NmutE) points]
\label{prop:NmutE}
An equilibrium point with $I_0^* > 0, I_1^*=0$, that is,
with no mutant-infected individuals and survival master-infected individuals, exists if and only if the master induced mortality vanishes $\delta_0=0$. 

If so, the corresponding microscopic equilibria are of the form
\[
 \QSmicro^{(1)}: \ (g_0^*, g_1^* \,; v_0^*, v_1^*) = \left( 0, 1, 0, \frac{\xi_1}{\gamma_1} \right) ,
\]
or
\[
\QSmicromu: \ (g_0^*, g_1^* \,; v_0^*, v_1^*) = \left( 1-\mu, \mu, \frac{\xi_0}{\gamma_0}(1-\mu), \frac{\xi_1}{\gamma_1}\mu\right), \qquad 0 < \mu \leq 1.
\]
\\
and the equilibria of the macroscopic system must fall into
one of the following cases:
\begin{itemize}
\item[(i)] Case $\chi>0$, \textbf{waning immunity}. There are four different scenarios according to the value of the transmission rates $\beta_{00}$ and $\beta_{01}$:
\begin{itemize}
\item[(i$_1$)] If $\beta_{00}=0$ and $\beta_{01}>0$ then it follows that, necessarily, the master recovery rate must vanish, $\pi_0=0$; the macroscopic equilibrium point takes the form
\[
(S^*, I_0^*, I_1^*, R^*, D^*) = (0,I_0^*,0,0,1-I_0^*).
\]

\item[(i$_2$)] If $\beta_{00}>0$ and $\beta_{01}=0$, then the macroscopic equilibrium reads
\[
(S^*, I_0^*, I_1^*, R^*, D^*) =
\left( \frac{\pi_0}{\beta_{00}}, I_0^*, 0, \frac{\pi_0}{\chi} I_0^*, 1-\frac{\pi_0}{\beta_{00}} - I_0^* \left( 1 + \frac{\pi_0}{\chi} \right) \right),
\]
with $0<I_0^*\leq \dfrac{1-\pi_0/\beta_{00}}{1+\pi_0/\chi}$.

\item[(i$_3$)] If $\beta_{00}>0$ and $\beta_{01}>0$, it follows that $\pi_0=0$ and the macroscopic equilibrium is
\[
(S^*, I_0^*, I_1^*, R^*, D^*) = \left( 0,I_0^*,0,0,1-I_0^* \right).
\]

\item[(i$_4$)] If $\beta_{00}=0$ and $\beta_{01}=0$, then, again, necessarily $\pi_0=0$; the macroscopic equilibrium has the expression
\[
(S, I_0^*, I_1^*, R, D) = (S,I_0^*,0,0,D), \qquad \textrm{with arbitrary $S$ and $D$.}
\]
\end{itemize}
\medskip

\item[(ii)] Case $\chi=0$, \textbf{permanent immunity}: it necessarily implies the master recovery rate to vanish, $\pi_0=0$. 
This, in its turn, it is divided into four cases, since their microscopic states at equilibrium might differ. They depend on the transmission rates $\beta_{00}$ and $\beta_{01}$. Indeed:
\begin{itemize}
\item[(ii$_1$)] If $\beta_{00}=\beta_{01}=0$, then 
\[
(S, I_0^*, I_1^*, R, D) = (S,I_0^*,0,R,D), 
\qquad \textrm{with arbitrary $S, R$ and $D$.}
\]
\item[(ii$_2$)] Otherwise, they take the form
\[
(S^*, I_0^*, I_1^*, R, D) = (0,I_0^*,0,R,D), 
\qquad \textrm{with arbitrary $R$ and $D$.}
\]
\end{itemize}

\end{itemize}
The corresponding macroscopic and microscopic steady state solutions are linked through compatibility conditions imposed between both systems. In our particular case, these constraints depend on the value of $\beta_{00}$ and $\beta_{01}$ as follows:
\begin{itemize}
\item[(a)] If $\beta_{00}=0$ then either $a_0=0$ or $v_0^*=0$. The latter may occur either because $\xi_0=0$ or $g_0^*=0$, which in turn implies $\mu=1$ in the microscopic equilibrium $\QSmicromu$.
 
\item[(b)] If $\beta_{01}=0$ then either $a_1=0$ or $v_1^*=0$. The second case implies either $\xi_1=0$ or $g_1^*=0$, which is only possible in the extreme case $\mu=0$ (see Sec.~\ref{sec:burnoutVirus}).

\item[(c)] The case $\beta_{00}>0$ never applies, since it necessarily requires that $v_0^*>0$, which is incompatible with the microscopic state $\QSmicro^{(1)}$.
\end{itemize}
\end{prop}

\begin{remark}
Like in the latter proposition, any full equilibrium state emerges from the coupling between the macroscopic and microscopic equilibria. 
They are subjected to mutual constraints, each one imposing compatibility conditions on the other. The number of admissible configurations can be substantially reduced by fixing certain parameters in advance.     
\end{remark}

\bigskip

And last but not least, the type of equilibrium states in which both strains, master and mutant,  co-circulate. As mentioned in the previous remark, in most of their cases, some micro--macro compatibility conditions have to hold in order to be real equilibrium points.

\begin{prop}[Co-circulating strains equilibrium (CSE) points]
\label{prop:coex:eqpoints}
Let consider a co-circulating strains equilibrium point, \emph{i.e.} $I_0^*>0$ and $I_1^*>0$. Then, necessarily, both mortality rates must vanish, $\delta_0=\delta_1=0$.
If so, the corresponding equilibria for the microscopic system are of the form:
\[
\begin{array}{lll}
\QSmicromuc: & \ (g_0^*, g_1^* \, ; v_0^*, v_1^*)  =
\left(  1 - \dfrac{\mu}{\mucrit^*} , \dfrac{\mu}{\mucrit^*} \, ;
\dfrac{\xi_0}{\gamma_0}\, g_0^*, \dfrac{\xi_1}{\gamma_1}\, g_1^*
\right) &\qquad \textrm{if $0 < \mu<\mucrit^*$}  \\[2.5ex]
 & \ (g_0^*, g_1^* \, ; v_0^*, v_1^*)  =
\left( 0, 1 \, ; 0, \dfrac{\xi_1}{\gamma_1} \right) &\qquad \textrm{if $\mucrit^* \leq \mu \leq 1$},
\end{array}
\]
where $\mucrit^*$ is given in~\eqref{def:mucrit_star}.
Moreover, it follows that the transmission rates $\beta_{01}$ and $\beta_{11}$ either vanish simultaneously or are both nonzero. The complete macroscopic and microscopic equilibria must fall into one of the following cases:
\begin{itemize}
\item[(i)] Case $\chi>0$, \textbf{waning immunity},
with two scenarios according to the value of the transmission rate $\beta_{00}$: 
\begin{itemize}
\item[(i$_1$)] If $\beta_{00}=0$ then, necessarily, the master recovery rate must vanish, $\pi_0=0$. Moreover we have:
\begin{itemize}
\item[(a)] If $\beta_{01}>0, \beta_{11}>0$ they are of the form
\[
(S^*,I_0^*, I_1^*, R^*, D)   \times (g_0^*, g_1^*; v_0^*, v_1^*) =
\left( \frac{\pi_1 I_1^*}{\beta_{01}^*I_0^* + \beta_{11}^*I_1^*} , I_0^*, I_1^*, \frac{\pi_1}{\chi} I_1^*, D^* \right) \times \QSmicromuc.
\]

\item[(b)] If $\beta_{01}=\beta_{11}=0$ then
there must be no mutant recovery $\pi_1=0$, and they are
\[
(S,I_0^*, I_1^*, R, D)   \times (g_0^*, g_1^*; v_0^*, v_1^*) = 
(S, I_0^*, I_1^*, 0, D) \times \QSmicromuc,
\]
with arbitrary $S$, $D$, and provided $\QSmicromuc$ satisfies $\beta_{01}^*=\beta_{11}^*=0$. 
\end{itemize}

\medskip

\item[(i$_2$)] If $\beta_{00}>0$ and the condition
\[
\pi_1 > \frac{\beta_{11}}{\beta_{00}} \pi_0
\]
holds, then
it follows that $\beta_{01}^*>0$ and $\beta_{11}^*>0$, and
the equilibrium points take the form
\begin{eqnarray*}
&&{(S,I_0^*, I_1^*, R, D) \times (g_0^*, g_1^*; v_0^*, v_1^*) = }\\[1.3ex] 
&&\qquad  \left( \frac{\pi_0}{\beta_{00}^*}, I_0^*, \frac{\beta_{01}^*\pi_0}{\beta_{00}^*\pi_1-\beta_{11}^*\pi_0} \, I_0^*, \frac{1}{\chi}\left(\pi_0+\pi_1\frac{\beta_{01}^*\pi_0}{\beta_{00}^*\pi_1-\beta_{11}^*\pi_0}\right) \, I_0^*, D^*\right) \times \QSmicromuc
\end{eqnarray*}
\end{itemize}

\medskip 

\item[(ii)] Case $\chi=0$, \textbf{permanent immunity}. 
Necessarily, both strains must have null recovery rate, that is $\pi_0=\pi_1=0$. Then, several scenarios arise:
\begin{itemize}
\item[(ii$_1$)] If $\beta_{00}=0$ we have two possible situations:
\begin{itemize}
\item[(a)] If $\beta_{01}>0, \beta_{11}>0$ then the equilibria are of type:
\[
(S^*,I_0^*, I_1^*, R, D)   \times (g_0^*, g_1^*; v_0^*, v_1^*) = 
(0, I_0^*, I_1^*, R, D) \times \QSmicromuc,
\]
with arbitrary $R$, $D$, and provided that $\QSmicromuc$ is compatible with the condition $\beta_{00}=0$.
\item[(b)] If $\beta_{01}=\beta_{11}=0$ they take the form
\[
(S,I_0^*, I_1^*, R, D)   \times (g_0^*, g_1^*; v_0^*, v_1^*) = 
(S, I_0^*, I_1^*, R, D) \times \QSmicromuc,
\]
with arbitrary $S$, $R$, and $D$, and provided that $\QSmicromuc$ is compatible with the conditions $\beta_{01}=\beta_{11}=0$.
\end{itemize}

\item[(ii$_2$)] If $\beta_{00}>0$, then
\[
(S^*,I_0^*, I_1^*, R, D)   \times (g_0^*, g_1^*; v_0^*, v_1^*) = 
(0, I_0^*, I_1^*, R, D) \times \QSmicromuc,
\]
with arbitrary $R$, and $D$.
\end{itemize}
\end{itemize}
\end{prop}

\begin{figure}[!ht]
    \centering
    \includegraphics[width = \linewidth]{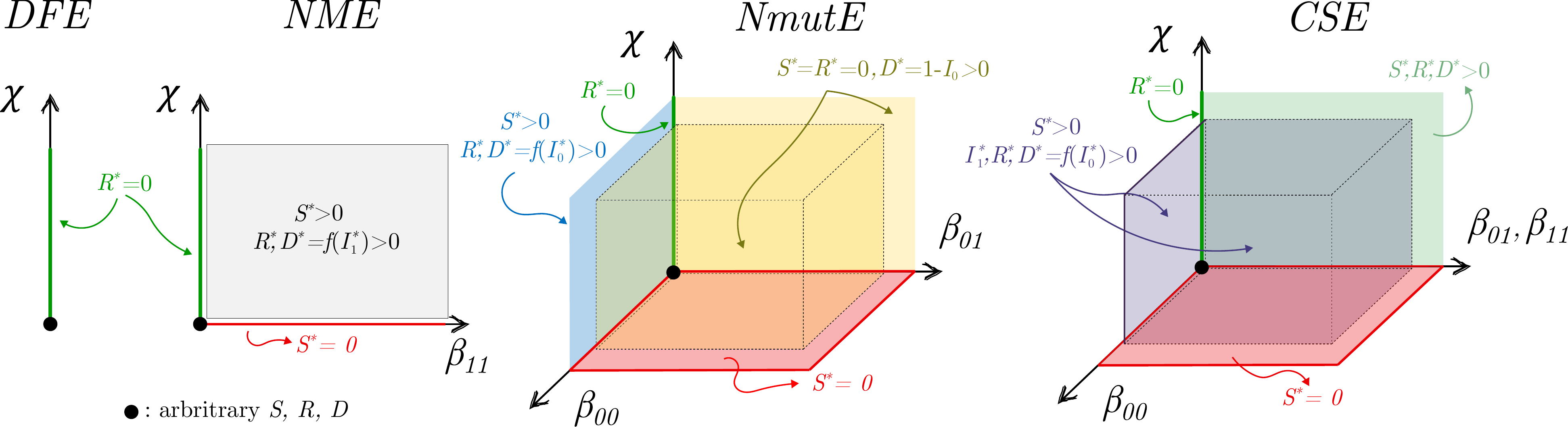}
    \caption{Schematic representation of the four different equilibrium points (DFE, NME, NmutE, and CSE) in terms of the parameters $\beta_{00}, \beta_{01}, \beta_{11}$, and $\chi$. Dashed lines represent open faces of the prism not included in the corresponding domain. The black dot indicates the origin.}
    \label{fig:EQS_summary}
\end{figure}

\begin{remark}
It is easy to see, from equation $\dot{D}=\delta_0 I_0 + \delta_1 I_1$, that the only possible scenario admitting - if any - periodic solutions is when both strains (master and mutant) do not induce mortality. Indeed, $\dot{D}>0$ and so $D$ increases unless $\delta_0=\delta_1=0$. This is clearly an interesting point to study but it is out of the scope of this work.
\end{remark}
Figure~\ref{fig:EQS_summary} provides a schematic representation of the system's equilibrium structure as a function of the transmission rates $\beta_{00}$, $\beta_{01}$,  $\beta_{11}$, and the waning immunity parameter $\chi$. Lines, planes and prisms are labeled with the corresponding equilibrium point form they give rise to. Dashed lines mean ``open" in the space of parameters (\emph{i.e.}, they do not include the boundary). 


\subsection[Computation of the basic reproduction number R0]{Computation of the basic reproduction number $\mathcal{R}_0$}
\label{sec:basic:reproduction:number}

The basic reproduction number $\mathcal{R}_0$ is a threshold quantity that measures the potential for the spread of an infectious disease. It is usually defined as the average number of secondary infections generated by a single infectious individual introduced into a fully susceptible population (see \cite{Hethcote00} and references therein). The computation of $\mathcal{R}_0$ provides an equivalent criterion, in terms of the parameters of the model, for determining the local stability of the DFE prior to the introduction of the first infectious individual. In our model~\eqref{eq:sirs:S}–\eqref{eq:sirs:D}, the DFE corresponds to the state
\[
(S,I_0,I_1,R, D)=(S(0),0,0,0,0)=(1,0,0,0,0),
\]
in which the population consists entirely of susceptible individuals (with no recovered or deceased individuals, as the disease has not yet emerged).

To compute $\mathcal{R}_0$ we use the so-called Next Generation Matrix (NGM) method, introduced in~\cite{DiekmannHeesterbeek00}. The NGM method we use in this work is based on~\cite{DriesscheWatmough02,DriesscheWatmough08}.
To this end, we split the variables $(S,I_0,I_1,R,D)$ into two compartments: the disease compartment $x=(I_0,I_1)$, and the disease-free compartment $y=(S,R,D)$. We then write the original system in the form
\begin{equation}
\dot{x}_i=\mathcal{F}_i(x,y) - \mathcal{V}_i(x,y) \quad i=1,\ldots,n, \qquad \qquad \dot{y}_j=g_j(x,y), \quad j=1,\ldots,m,
\label{NGM:subsystems:DW}
\end{equation}
where $\mathcal{F}_i$ denotes the inflow of new infections in compartment $i$, and $\mathcal{V}_i$ the rate of transmissions between compartment $i$ and other infected compartments. As noted explicitly in~\cite{DriesscheWatmough08}, the decomposition into $\mathcal{F}$ and $\mathcal{V}$ is not unique and depends on the biological interpretation of the compartments and on whether the transitions correspond to infection events or transmissions among infected states. In this notation, the DFE is $(0,y_0)$. 

The required NGM assumptions are checked in Appendix~\ref{app:NGMdetails}. Hence, the linearization of the infected subsystem at the DFE $(0,y_0)$ reads
\begin{equation}
\dot{x}=(F-V)x, \qquad \textrm{where} \quad 
F=\frac{\partial \mathcal{F}_i}{\partial x_j}(0,y_0), \qquad 
V=\frac{\partial \mathcal{V}_i}{\partial x_j}(0,y_0), \quad \textcolor{blue}{i,j=1,\ldots,n}.
\label{NGM:F:V:ode}
\end{equation}
It was shown in~\cite{DiekmannHeesterbeek00} that the local stability of system $\dot{x}=\mathcal{F}(x,y) - \mathcal{V}(x,y)$ around the DFE $(0,y_0)$ is determined by the linear stability of~\eqref{NGM:F:V:ode}. Furthermore, the matrix $FV^{-1}$ is referred to as the next generation matrix. Defining $\mathcal{R}_0=\rho(FV^{-1})$, as its spectral radius, the following result holds:
the DFE $(0, y_0)$ is locally asymptotically stable if $\mathcal{R}_0 < 1$ and unstable if $\mathcal{R}_0 > 1$.
In our case, the disease compartment satisfies the ODE
\begin{equation}
\left(
\begin{array}{c}
\dot{I}_0 \\ \dot{I}_1    
\end{array}
\right) = \left(
\begin{array}{c}
\beta_{00} SI_0 \\ \beta_{01} S I_0 + \beta_{11} S I_1
\end{array}
\right) -
\left(
\begin{array}{c}
(\pi_0+\delta_0) I_0 \\ (\pi_1+\delta_1) I_1
\end{array}
\right) = \mathcal{F}(X) - \mathcal{V}(X),
\label{NGM:disease:comp:ode}    
\end{equation}
and the disease-free subsystem is given by
\begin{equation}
\left(
\begin{array}{c}
\dot{S} \\ \dot{R} \\ \dot{D}
\end{array}
\right) = \left( 
\begin{array}{c}
- (\beta_{00}+\beta_{01}) S I_0 - \beta_{11} SI_1 + \chi R \\
\pi_0 I_0 + \pi_1 I_1 - \chi R \\
\delta_0 I_0 + \delta_1 I_1
\end{array}
\right) = g(X),
\label{NGM:disease-free:comp:ode}    
\end{equation}
where $X=(S,I_0,I_1,R,D)$.
Consequently, the NGM is
\[
K=F V^{-1} = \left( 
\begin{array}{cc}
\dfrac{\beta_{00}}{\pi_0 + \delta_0} & 0 \\[1.4ex]
\dfrac{\beta_{01}}{\pi_0 + \delta_0} & \dfrac{\beta_{11}}{\pi_1 + \delta_1}
\end{array}
\right),
\]
and the basic reproduction number is 
\[
\mathcal{R}_0 = \rho(FV^{-1}) = \max \left\{ \frac{\beta_{00}}{\pi_0 + \delta_0}, 
\frac{\beta_{11}}{\pi_1 + \delta_1} \right\},
\]
which corresponds to the maximum between the basic reproduction numbers of both viral strains when they are assumed to occur independently. However, recall that in the DFE case, $I_0=I_1=0$, we have defined $\nu_0=0$ and $\nu_1=0$ and therefore
\[
\beta_{01}=\beta_{01}(v_1(0))=\frac{a_1 \nu_0 v_1(0)}{b_1 + \nu v_1(0)}=0 \quad  \textrm{and} \quad \beta_{11}=\beta_{11}(v_1(0))=\frac{a_1 \nu_1 v_1(0)}{b_1 + (1-\nu)v_1(0)}=0.
\]
Hence, the basic reproduction number is given by
\begin{equation}
\mathcal{R}_0 = \rho(FV^{-1}) = \frac{\beta_{00}}{\pi_0 + \delta_0},
\label{def:R0}
\end{equation}
the same as the one corresponding to the master strain. This fact agrees with the intuitive idea that at the very starting point of the epidemics, the role played by the mutant strain is, for a while, secondary. The quotient
$\beta_{00}/(\pi_0 + \delta_0)$
represents the proportion between the infectious $I_0$-inflow ($\beta_{00}$, infection) and outflow ($\pi_0+\delta_0$ immunization and death) rates. However, this criterion only describes the initial stage, when the viral population is effectively monomorphic. As the epidemic evolves, mutation and transmission generate a mixture of strains, and the dynamics can no longer be captured by the master strain alone. This motivates the study of the condition for a pandemic growth that extends the inflow--outflow balance to account for the combined contribution of both strains and their cross-scale interaction.


\subsection{On the condition for pandemic growth}
\label{se:growing_pandemic}

In analogy with epidemic models based on the instantaneous effective reproduction number - often denoted by $\mathcal{R}_t$ when it varies in time - the goal of this section is to derive a condition that ensures pandemic growth and to relate it to the dynamics of both infected populations. Focusing first on individuals infected by the master strain, we obtain
\[
\dot{I}_0=\left( \beta_{00} S - (\pi_0 +\delta_0) \right) I_0,
\]
which leads to the definition of
\begin{equation*}
\mathcal{R}_t^{(0)}= \frac{\beta_{00}}{\pi_0 + \delta_0} S,   \qquad \textrm{if} \ I_0 \ne 0,
\label{def:Rt:0}
\end{equation*}
where $S=S(t)$ and 
\[
\beta_{00}=\frac{a_0 v_0(t)}{b_0 + v_0(t)}.
\]

This is the standard definition for the (time) reproduction number in a one-strain epidemic. It satisfies that $I_0$ grows if and only if $I_0\ne 0$ and $\mathcal{R}_t^{(0)}>1.$ We can analogously seek for a condition of $I_1$-infected individuals. Thus, for $I_1\ne 0$, it follows that
\[
\dot{I}_1= \left( \beta_{01} I_0 + \beta_{11} I_1 \right) S - (\pi_1 + \delta_1) I_1 > 0
\]
is equivalent to $\mathcal{R}_t^{(1)}>1$, for
\begin{equation}
\mathcal{R}_t^{(1)} = \frac{1}{\pi_1 + \delta_1} \left( \beta_{01} \frac{I_0}{I_1} + \beta_{11} \right) S,
\label{def:Rt:1}    
\end{equation}
where $S=S(t)$, $I_j=I_j(t)$, $j=0,1$, and
\[
\beta_{01}=\frac{a_1 \nu v_1(t)}{b_1 + \nu v_1(t)} \quad \mathrm{and} \quad \beta_{11}=\frac{a_1 (1-\nu) v_1(t)}{b_1 + (1-\nu) v_1(t)},
\]
where recall that we denote $\nu_0=\nu$ and $\nu_1=1-\nu$ when $0<\nu_0<1$.
Let us now compare the conditions above with the one ensuring the growth of the pandemic whatever the strain (master or mutant) is causing it. This is equivalent to impose that $\dot{I}_0 + \dot{I}_1 >0$. Indeed,
\begin{eqnarray}
\dot{I}_0 + \dot{I}_1 &=& \bigg( (\beta_{00}+\beta_{01}) I_0 + \beta_{11} I_1 \bigg) S - (\pi_0 + \delta_0) I_0 - (\pi_1 + \delta_1) I_1 > 0 \nonumber \\
&\Leftrightarrow &  \mathcal{G}_t=(\pi_0 + \delta_0) I_0 \left( \mathcal{R}_t^{(0)} -1 \right) + (\pi_1 + \delta_1) I_1 \left( \mathcal{R}_t^{(1)} - 1 \right) > 0.
\label{def:growing:I:condition}
\end{eqnarray}
Notice that the terms $(\pi_j+\delta_j)I_j$, for $j=0,1$, in the expression~\eqref{def:growing:I:condition} above correspond to the outflow (either by immunisation or death) of infectious population density $I_j$.
If we denote by
\begin{equation*}
\mathcal{O}_t^{(j)} = \left( \pi_j + \delta_j \right) I_j(t), \qquad j=0,1,    
\end{equation*}
the outflow rate at time $t$ of $I_j$, then $\mathcal{G}_t>0$ is equivalent to 
\begin{equation}
\mathcal{G}_t = \mathcal{O}_t^{(0)} \left( \mathcal{R}_t^{(0)} -1 \right) + \mathcal{O}_t^{(1)} \left( \mathcal{R}_t^{(1)} -1 \right) >0. 
\label{def:G_t}
\end{equation}

\begin{remark}
Let us denote by
\[
R_t^{(1)} = \frac{\beta_{11}}{\pi_1+\delta_1} S(t)
\]
the corresponding $R_t$ value of a pandemic driven uniquely by the mutant variant ($v_1$, and so $I_1$ only). From~\eqref{def:Rt:1} and the expression above, it turns out that
\[
\mathcal{R}_t^{(1)} = \frac{\beta_{01}}{\pi_1 + \delta_1} \frac{I_0}{I_1} + R_t^{(1)},
\]
which relates the $R_t$ values of a $v_1$ mutant driven pandemic when $I_1$ acts alone, \emph{i.e.} $R_t^{(1)}$, and when it comes from a mutation/competition relation with its master virus sequence, \emph{i.e.} $\mathcal{R}_t^{(1)}.$ Notice that 
\[
\frac{\mathcal{R}_t^{(1)}}{R_t^{(1)}} = 1 + \frac{\beta_{01}}{\beta_{11}} \frac{I_0}{I_1} > 1
\]
being, therefore 
\[
\mathcal{R}_t^{(1)} > R_t^{(1)},
\]
both of them computed assuming the same quantity $S(t)$ of susceptible individuals at time $t$.
\end{remark}

\subsection{Case of interest 1: a vaccine-like viral strain.}
\label{sec:COI1}

One illustrative scenario chosen to demonstrate the applicability of the model concerns the emergence of a mutant strain that elicits a faster recovery than its master strain, \emph{i.e.} $\pi_1 > \pi_0$, while both strains are nonlethal ($\delta_0 = \delta_1 = 0$). Several real examples fall within this vaccine‑like regime. The best-documented case is Sabine's live attenuated oral poliovirus vaccine type 2 (OPV2), which replicates in the gut, induces strong mucosal immunity, and can secondarily spread to close contacts, thereby boosting herd immunity~\cite{Macklin2023,Gast2022}. This example maps directly onto our case assumptions: low $\delta_1$, high $\pi_1$, and transmission rates tied to fecal viral load, such that larger $v_1$ drives $\beta_{11}$ toward saturation. OPV’s historic impact stems from this benign circulation, tempered by the rare risk of genetic reversion to circulating vaccine-derived poliovirus (cVDPV) —exactly the mutation-selection tension our framework captures via $\mu_c(\nu)$. From 2021–2025, Novel OPV2 (nOPV2) was deployed at scale to preserve OPV’s mucosal and transmission benefits while improving genetic stability. Field and laboratory data show substantially fewer emergences than OPV2, consistent with an avirulent, strongly immunizing agent that can circulate without appreciable mortality, squarely “vaccine‑like” \cite{Estivariz2023}.

A second example is the live‑attenuated rotavirus vaccines. They are shed in stool following the first dose, and although documented transmission to contacts is uncommon, it can confer indirect protection \cite{Zhou2019}. Recent studies link shedding to seroconversion, and surveillance data report very low transmission without symptomatic disease in exposed infants \cite{Zalot2025}, again matching high $\pi_1$, $\delta_1 \approx 0$, and a typically low‑to‑moderate, short‑lived $\beta_{11}(v_1)$ as the recovered population accumulates.

The assumptions $\pi_1 > \pi_0$ and $\delta_0 = \delta_1 = 0$ allow, in principle (Section~\ref{sec:eqPoints}), the existence of all four types of steady states: DFE, NME, NmutE, and CSE equilibria. Another direct implication of $\delta_0 = \delta_1 = 0$ is that $\dot{D}=\delta_0 I_0 + \delta_1 I_1 =0$ and so $D$ is a first integral. Recall that the total mass $S+I_0+I_1+R+D$ is also a first integral of the macroscopic system. Since it is expected that, at the beginning of the process, $D(0)=0$, henceforth it will be assumed that $D(t)=0$ $\forall t\geq 0$ and $D$ will be excluded from further analysis.

\begin{figure}[!ht]
    \centering
    \includegraphics[width = \textwidth]{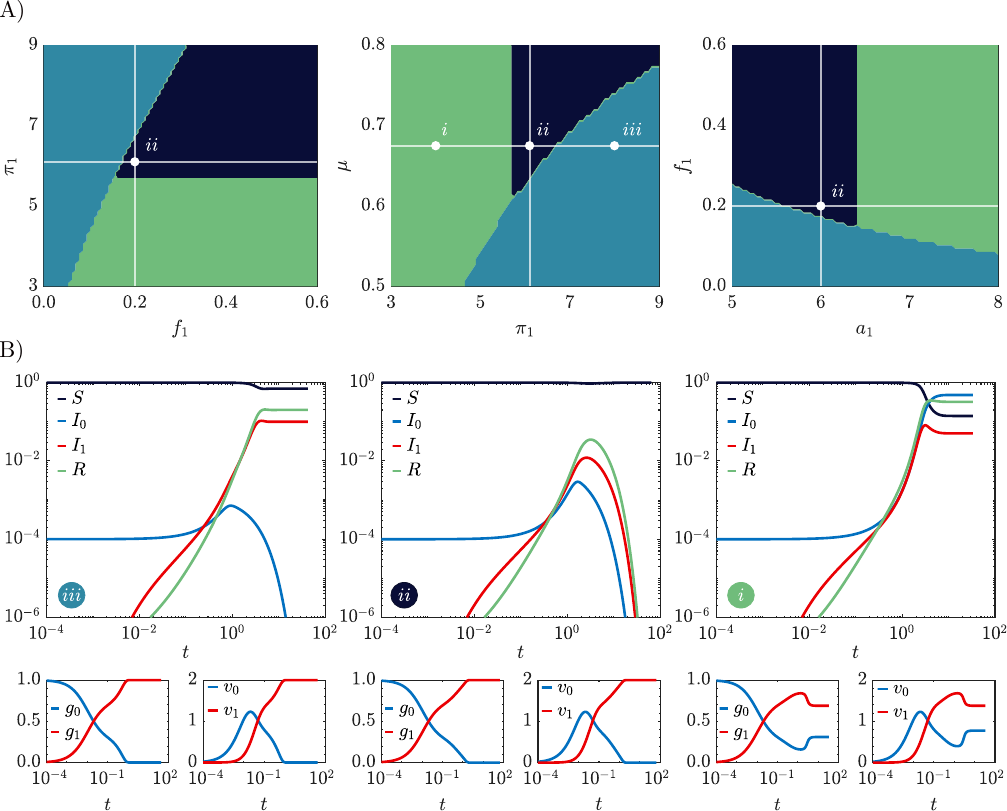}
    \caption{(A) $\omega$-limits of the principal trajectory with i.c.~\eqref{case:interest:2:ic} and fixed parameters~\eqref{caseinterest1:CSE:parameters} for three different combinations of parameters. Each point of the plot is colored according to the type of equilibrium reached: \textcolor{ColorNME}{NME}, \textcolor{ColorDFE}{DFE}, and \textcolor{ColorCSE}{CSE}. The first and second panels correspond to $a_1 = 6$, while the third one fixes $\mu = 0.675$. These, together with $f_1 = 0.2$, will be the nominal values employed for the case study. They are highlighted in the plots using solid white lines to facilitate visual reference. (B) Time evolution of the macroscopic and microscopic variables for three scenrios, highlighted in the upper panel. Each simulation illustrates convergence toward a different equilibrium.}
    \label{fig:CoI3_highIRmut}
\end{figure}

For a first numerical exploration, we have conveniently fixed the following parameters:
\begin{eqnarray}
&&\chi = 2,\quad \pi_0 = 0.5,\quad  f_0 = 1,\quad  \xi_0 = 2,\quad  \gamma_0 = 0.8, \label{caseinterest1:CSE:parameters} \\
&&  a_0 = 4,\quad  b_0=0.1, \quad \xi_1 = 1,\quad  \gamma_1 = 0.5,\quad b_1 = 0.1, \nonumber
\end{eqnarray}
with $\delta_0=\delta_1=0$ and the ratio between the slow and fast time scales equal to $\eps=0.01$. For values of the recovery rate $\pi_1>\pi_0$, the mutant fitness $f_1$, the maximal potential threshold for $\beta_{01}$ and $\beta_{11}$, and the master's mutation error $\mu$, ranging in suitable intervals, we compute the $\omega$-limit of the trajectory starting with initial conditions
\begin{equation}
(S,I_0,I_1,R,D) \times (g_0,g_1; v_0, v_1) = (1-10^{-4}, 10^{-4}, 0,0,0) \times (1,0;0,0).
\label{case:interest:2:ic}
\end{equation}
Along the paper, we will refer to the solution with these i.c. as the principal trajectory. According to the type of equilibrium point reached (its $\omega$-limit), for any choice of the latter parameters, we paint (Figure~\ref{fig:CoI3_highIRmut}A) the corresponding point in the parameter's space with colors:
\textcolor{ColorNME}{NME}, \textcolor{ColorDFE}{DFE}, and \textcolor{ColorCSE}{CSE}. In Figure~\ref{fig:CoI3_highIRmut}B the time evolution of the macroscopic and microscopic variables, respectively, for three particular choices of these parameters is shown. 
To analyze the variation of the $\omega$-limit (and its time evolution) of this principal trajectory in terms of the mutant recovery $\pi_1$, we fix the rest of parameters as $a_1 = 6$, $f_1 = 0.2$, and $\mu = 0.675$.

Figure~\ref{fig:CoI3_highIRmut} suggests the following points to investigate in terms of $\pi_1$: (a) the existence of equilibrium points and their types; (b) their local stability and multistability scenarios and, when possible, (c) to find out mechanisms underlying changes in the $\omega$-limit of the principal trajectory with i.c.~\eqref{case:interest:2:ic}.

\subsubsection{Equilibrium points}
The existence and expression of these equilibria follow from the propositions in Section~\ref{sec:eqPoints}. Further analysis will show how variations in the host immune response (encoded in the value of $\pi_1$) influence their behavior. We will particularly be concerned with the possibility of co-circulation of both the master and the mutant viral strains. 
Recall that $D\equiv 0$ will be always assumed.
Thus, regarding the equilibrium points in this particular case, it follows that:

\noindent\textbf{$\bullet$ Disease-free equilibrium (DFE)} \\[1.1ex]
They are of the form
\begin{equation*}
\textrm{DFE}: \quad (S^*,I_0^*,I_1^*,R^*,D^*)\times(g_0^*,g_1^*,v_0^*,v_1^*)=\left(1,0,0,0,0\right)\times\left(g_0^*,1-g_0^*,\frac{\xi_0}{\gamma_0}g_0^*,\frac{\xi_1}{\gamma_1}(1-g_0^*)\right), 
\end{equation*}
for $g_0^* \in [0,1]$. However, they represent a segment in the full system space, and so there are infinitely many of them, in terms of the macroscopic equilibrium they correspond to the same steady state $(S^*,I_0^*,I_1^*,R^*,D^*)=(1,0,0,0)$.

\noindent\textbf{$\bullet$ No master equilibrium (NME)}. \\[1.1ex]
From Proposition~\ref{prop:NME}, the only microscopic equilibrium is given by
$(g_0,g_1)=(0,1)$ which, in its turn, yields to
\[
(v_0,v_1) = \left( 0, \frac{\xi_1}{\gamma_1} \right) \qquad \textrm{and} \qquad \beta_{00} =0, \quad  \beta_{01}=0, \quad \beta_{11} = \frac{a_1 \frac{\xi_1}{\gamma_1}}{b_1 + \frac{\xi_1}{\gamma_1}}.
\]
From the expression for the $I_1$ equilibrium, the following constraint is derived:
\[
I_1^* = \frac{1- \frac{\pi_1}{\beta_{11}}}{1+\frac{\pi_1}{\chi}} \in (0,1)
 \Rightarrow 
 1-\frac{\pi_1}{\beta_{11}} > 0 \Rightarrow \pi_1 < \beta_{11} \Leftrightarrow \pi_1 < \frac{a_1 \xi_1}{b_1\gamma_1 + \xi_1},
\]
which in our particular case reads $\pi_1<40/7$.
If such condition holds, the expression for the unique NME equilibrium is
\begin{equation*}
\textrm{NME}:\quad(S^*,I_0^*,I_1^*,R^*,D^*)\times(g_0^*,g_1^*; v_0^*,v_1^*)=
\left(
\frac{\pi_1}{\beta_{11}}, 0 , 
\frac{1- \frac{\pi_1}{\beta_{11}}}{1+\frac{\pi_1}{\chi}},
\frac{\pi_1}{\chi} \cdot \frac{1- \frac{\pi_1}{\beta_{11}}}{1+\frac{\pi_1}{\chi}},
 0 \right)  \times   \left(0,1;0,\frac{\xi_1}{\gamma_1}\right).
\end{equation*}

\noindent\textbf{$\bullet$ No mutant equilibrium (NmutE)} \\[1.1ex]
There are no equilibrium points of this type. Briefly, the only possible microscopic solutions $(g_0,g_1)=(1-\mu,\mu)$ with $0\leq \mu <1$ lead to
\[
(v_0,v_1) = \left( \frac{\xi_0}{\gamma_0}(1-\mu), \frac{\xi_1}{\gamma_1} \mu \right) \qquad  \textrm{and} \qquad 
\beta_{00} = \frac{a_0 \frac{\xi_0}{\gamma_0}(1-\mu)}{b_0 + \frac{\xi_0}{\gamma_0}(1-\mu)}, \qquad
\beta_{01} = \frac{a_1 \frac{\xi_1}{\gamma_1} \mu}{b_1 + \frac{\xi_1}{\gamma_1}\mu}, \qquad \beta_{11}=0.
\]
From the third equation of the macroscopic system, $\beta_{01}I_0 S=0$, it turns out that either $S=0$ or $\beta_{01}=0$ (\textit{i.e.}, $a_1=0$). The first of both cases, $S=0$, does not lead to any equilibrium solution since it leads to $\pi_0 I_0=0$, which is not possible. On the other hand, from the third equation of the macroscopic system, $\beta_{01}I_0 S=0$, it follows that either $S=0$ or $\beta_{01}=0$. The first case, $S=0$ implies again that $\pi_0 I_0=0$, a contradiction. The second case is only possible if $a_1=0$ (a non-transmissible mutant), which is not the case under consideration.

\noindent\textbf{$\bullet$ Co-circulating strains equilibrium (CSE)} \\[1.1ex]
This is the most general and intrincated case, as it is subject to the fewest restrictions. We seek for equilibrium points with both infected populations nonvanishing. This implies that
$\nu^* =\nu(I_0^*, I_1^*)\in (0,1)$ and, according to the Swetina-Schuster model, the genomes' equilibrium are given either by
\[
(g_0^*,g_1^*) = \left( 1 - \frac{\mu}{\mucrit^*}, \frac{\mu}{\mucrit^*} \right) \qquad \textrm{if $0<\mu < \mucrit^*$},
\]
or $(g_0,g_1)=(0,1)$ (master genome's extinction) if $\mucrit^* \leq \mu \leq 1$, where we recall that
\begin{equation*}
\mucrit^* = 1 - \frac{f_1}{f_0} \left( \frac{1}{\nu^*} - 1 \right)
\end{equation*}
Let us analyze both cases separately.
\begin{itemize}
\item[(\emph{i})] Case $0< \mu < \mucrit^*$.\\[1.1ex]
Here we have:
\[
(g_0^*,g_1^*) = \left( 1 - \frac{\mu}{\mucrit^*}, \frac{\mu}{\mucrit^*} \right), \qquad \quad (v_0^*,v_1^*) = \left( \frac{\xi_0}{\gamma_0}g_0^*, \frac{\xi_1}{\gamma_1} g_1^* \right),
\]
and
\[
\beta_{00}^* = \frac{a_0 v_0^*}{b_0 + v_0^*}, \qquad
\beta_{01}^* = \frac{a_1 \nu^* v_1^*}{b_1 + \nu^*v_1^*}, \qquad
\beta_{11}^* = \frac{a_1 (1-\nu^*) v_1^*}{b_1 + (1-\nu^*)v_1^*},
\]
Notice that $\beta_{00}^*$, $\beta_{01}^*$, and $\beta_{11}^*$ are all three strictly positive.
On the other side, from Proposition~\ref{prop:coex:eqpoints} the following expressions for the macroscopic variables at equilibrium hold:
\begin{equation}
S^* = \frac{\pi_0}{\beta_{00}^*} \qquad  R^* = \frac{\pi_0 I_0^*+\pi_1 I_1^*}{\chi}, \qquad I_1^* = I_0^*\left(\dfrac{\beta_{01}^* \pi_0 }{\beta_{00} \pi_1^* -\beta_{11}^*\pi_0 }\right).
\label{caseInterest:1:CSE:equilibria}
\end{equation}
Hence, a necessary condition for the existence of $I_1^*$ (and, therefore, of CSE) is that
\begin{equation*}
\frac{\pi_0}{\beta_{00}^*} < \frac{\pi_1}{\beta_{11}^*}.
\label{CSE:condition:pi1}
\end{equation*}
If it holds, macroscopic CSE solutions must satisfy
\[
S^* + I_0^* + I_1^* + R^* = 1 \Leftrightarrow 
\frac{\pi_0}{\beta_{00}^*} + I_0^* + I_0^*\left(\frac{\beta_{01}^* \pi_0}{\beta_{00}^* \pi_1 -\beta_{11}^* \pi_0 }\right) +  \frac{\pi_0 I_0^*+\pi_1 I_1^*}{\chi} = 1
\]
or, equivalently, $\mathcal{F}(I_0^*, I_1^*)=0$, where 
\begin{equation*}
\mathcal{F}(I_0^*,I_1^*) =  \frac{\pi_0}{\beta_{00}^*} + I_0^* + I_0^*\left(\frac{\beta_{01}^*\pi_0}{\beta_{00}^*\pi_1 -\beta_{11}^*\pi_0 }\right) +  \frac{\pi_0 I_0^*+\pi_1 I_1^*}{\chi} - 1.   
\end{equation*}
This is an involved (but rational) equation in $I_0^*, I_1^*$ since $\beta_{01}^*$, and $\beta_{11}^*$, in turn, also depend on them through $\nu^*$. This equation has been numerically solved in Figure~\ref{fig:CEEpoints}).

\begin{figure}[!ht]
    \centering
    \includegraphics[scale = 1]{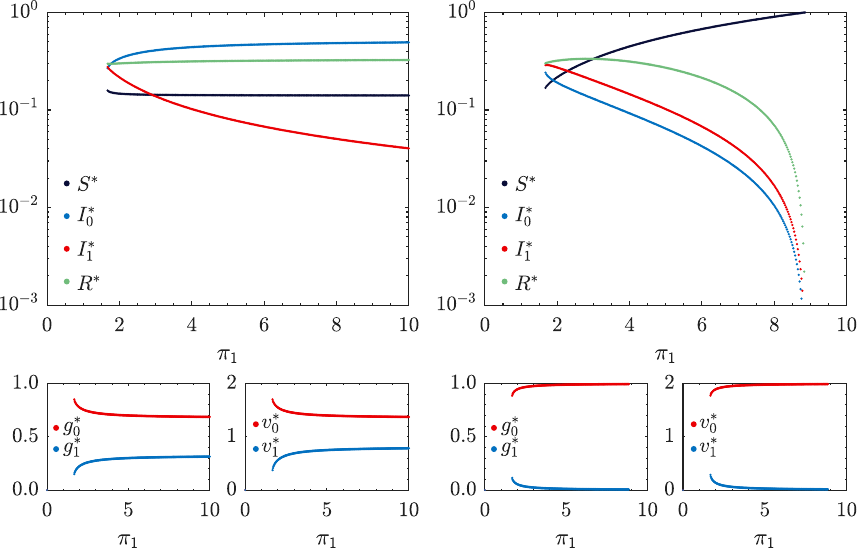}
    \caption{Co-circulating strains equilibrium points, computed numerically, as a function of $\pi_1\in (\pi_0,10)=(0.5,10)$. If we define $\pi_1^{\dagger}\simeq 1.675$ and $\pi_1^{\diamond}\simeq 8.875$ then: for 
    $\pi_1 \in (\pi_0,\pi_1^{\dagger})$ there is no equilibrium point; for $\pi_1\in (\pi_1^{\dagger},\pi_1^{\diamond})$ there are two, and for $\pi_1 \in (\pi_1^{\diamond},10)$ there is only one. They are plotted in the left and right panels above (top: macroscopic, bottom: microscopic). These two families of CSE equilibria show relevant differences: on one side, family CSE$_1$ (left) starts at a point with equal infective individuals $I_0^*=I_1^*$ and exhibits a monotonous decrease of the $I^*_1$-population as $\pi_1$ grows; on the other side, family CES$_2$ (right) begins with $I^*_1>I^*_0$ and undergoes a bifurcation (namely,  collision with the DFE equilibrium $(1,0,0,0)$) as $\pi_1 \rightarrow \pi_1^{\diamond}$ (see Figure~\ref{fig:placeholder}). 
    }  
    \label{fig:CEEpoints}
\end{figure}

\begin{remark}
At this CSE equilibrium, the $I_0$ and $I_1$-reproduction numbers (see Section~\ref{se:growing_pandemic}) at the equilibrium,
\[
\mathcal{R}_*^{(0)} = \frac{\beta_{00}^*}{\pi_0} S^*, \qquad 
\mathcal{R}_*^{(1)} = \frac{1}{\pi_1} \left( \beta_{01}^* \frac{I_0^*}{I_1^*} + \beta_{11}^* \right) S^*,
\]
are both equal to $1$. Indeed, first, having in mind~\eqref{caseInterest:1:CSE:equilibria}, it follows that
\[
\mathcal{R}_*^{(0)} = \frac{\beta_{00}^*}{\pi_0} S^* =  \frac{\beta_{00}^*}{\pi_0}  \frac{\pi_0}{\beta_{00}^*}=1.
\]
Regarding the second assertion:
\begin{equation*}
\mathcal{R}_*^{(1)} = \frac{\pi_0}{\beta_{00}^* \pi_1} \left( \beta_{01}^* \cdot \frac{\beta_{00}^* \pi_1 - \beta_{11}^*\pi_0}{\beta_{01}^* \pi_0} + \beta_{11}^* \right)
= \frac{1}{\beta_{00}^* \pi_1} \left( \beta_{00}^* \pi_1 - \beta_{11}^* \pi_0 + \beta_{11}^* \pi_0 \right) = 1.
\end{equation*}
\end{remark}
\medskip

\item[(\emph{ii})] Case $\mucrit \leq \mu \leq 1$. \\[1.1ex]
In this scenario, there are no equilibrium points of this type. Indeed, the microscopic equilibrium is $(g_0,g_1,v_0,v_1)=(0,1,0,\xi_1/\gamma_1)$,
which means $\beta_{00}=1$. But this implies $\dot{I}_0=-\pi_0I_0=0$ which can be only satisfied if $\pi_0=0$ or $I_0=0$, in contradiction with our assumptions.
\end{itemize}

\begin{figure}[!ht]
    \centering
    \includegraphics[scale=1]{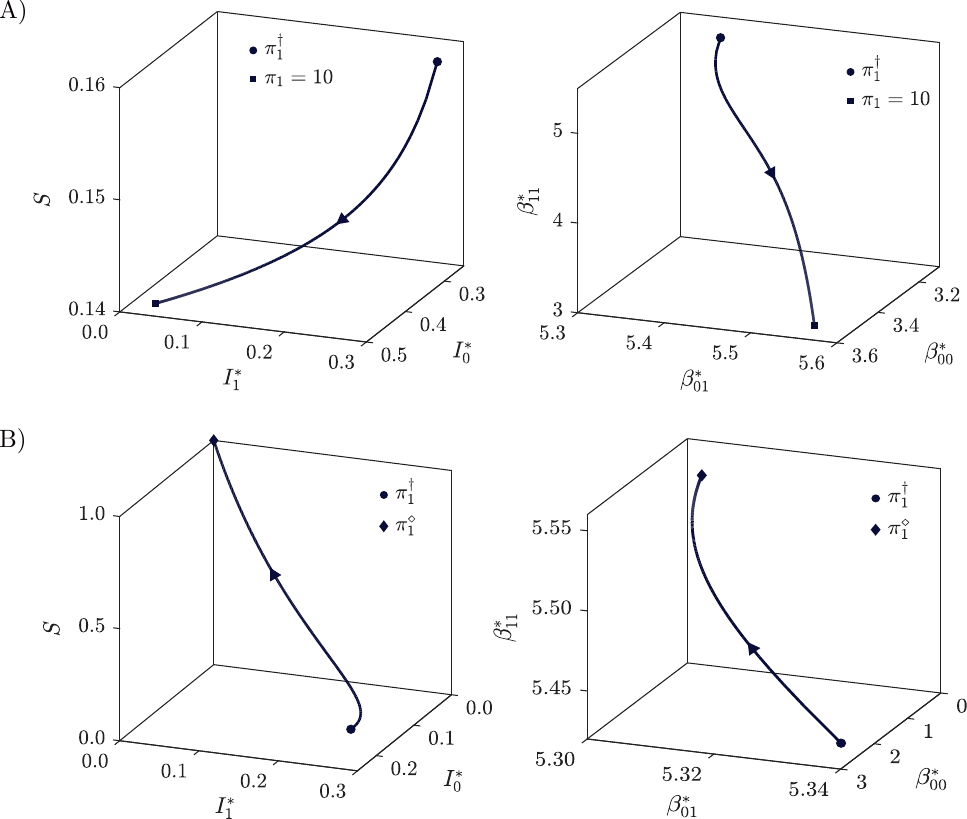}
    \caption{(A) Evolution of the family of CSE$_1$ points, varying with $\pi_1$ (see Figure~\ref{fig:CEEpoints} (left)), for the $(S,I_0,I_1)$ variables (left) and $\beta_{00}^*$, $\beta_{01}^*$ (right). (B) Mechanism leading the CSE$_2$ point (see Figure~\ref{fig:CEEpoints}(right)) towards its collision with the DSE point $(1,0,0,0)$ as $\pi_1$ tends to $\pi_1^{\diamond}$. Left: $(S,I_0^*,I_1^*)$ variables. Right: transmission rates $\beta_{00}^*$ and $\beta_{01}^*$.}
    \label{fig:placeholder}
\end{figure}

\subsubsection{Equilibria stability}
\label{se:caseofinterest1:equilibria}

In the previous section, we identified the admissible equilibrium states of the system under the parameter regime of interest. In this section, we analyze the local stability of each equilibrium point by evaluating the corresponding jacobian matrices. Recall that the parameter choices are specified in~\eqref{caseinterest1:CSE:parameters}, with $\delta_0=\delta_1=0$, $\eps=0.01$, and $\pi_1$ varying in the interval $(\pi_0,10).$

The Jacobian matrix determines the linearized dynamics of the system in the vicinity of an equilibrium, and its eigenvalues determine the local behavior. In particular, an equilibrium is locally asymptotically stable if and only if all eigenvalues of the jacobian have negative real parts. Conversely, the presence of any eigenvalue with a positive real part implies instability. This analysis allow us to identify parameter regions in which qualitative changes in stability occur, providing insights into potential bifurcations and transitions between different dynamical regimes.

\noindent\textbf{$\bullet$ Disease-free equilibria (DFE).}
They are of the form
\begin{equation*}
\textrm{DFE}: \quad (S^*,I_0^*,I_1^*,R^*,D^*)\times(g_0^*,g_1^*,v_0^*,v_1^*)=\left(1,0,0,0,0\right)\times\left(g_0^*,1-g_0^*,\frac{\xi_0}{\gamma_0}g_0^*,\frac{\xi_1}{\gamma_1}(1-g_0^*)\right), 
\end{equation*}
parametrized in terms of $g_0^*\in [0,1]$. Their jacobian have the following spectrum:
\begin{equation*}
    \lambda_\mathrm{DFE} = \bigg\{-2 ,\, -\frac{4}{5} ,\, -\frac{1}{2} ,\, 0 ,\, 0 ,\, 0 ,\, \frac{175g_0^*-1}{50g_0^*+2} ,\, -\pi_1\,\bigg\}.
\end{equation*}
All the eigenvalues are real. Three of them are $0$ (associated to the freedom in $I_0=I_1=D=0$), four negative and a unique eigenvalue
\[
\lambda_{\mathrm{dfe}}(g_0^*)=\frac{175g_0^*-1}{50g_0^*+2},
\]
exclusively depending on $g_0^*$, which determines their stability. This was expected since the macroscopic equilibrium is always $(S^*,I_0^*,I_1^*,R^*,D^*)=(1,0,0,0,0)$. The variation in the complete system equilibrium comes from the microscopic equilibria, 
given by a segment, parametrized by $g_0^*$, in the genomes-virions space. 
Namely, it is easy to check that an DFE is asymptotically stable (indeed, attractor in the normal bundle space to the invariant manifold $\{ I_0=I_1=D=0 \}$) provided that $g_0^*\in [0,1/175)$ and unstable (a saddle point) if $g_0^* \in (1/175,1]$. 

Notice that the stability condition is equivalent (as expected) to the basic reproduction number $\mathcal{R}_0<1$ defined in~\eqref{def:R0}. Indeed, from the fact that $\beta_{11}^*=0$, from the values of the parameters in~\eqref{caseinterest1:CSE:parameters}, and taking into account that $v_0^* =\xi_0g_0^*/\gamma_0$, it follows that
\[
\mathcal{R}_0  = \frac{\beta_{00}^*}{\pi_0} = \frac{20 g_0^*}{0.1 + 2.5 g_0^*.}
< 1 \Leftrightarrow g_0^* < \frac{1}{175},
\]
for any value of $\pi_1$.

\bigskip

\begin{figure}[!ht]
    \centering
    \includegraphics[scale = 1]{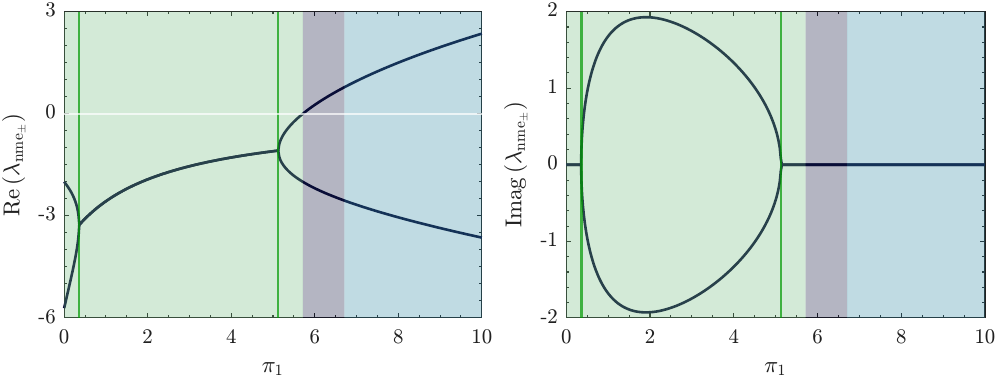}
    \caption{Eigenvalues spectrum of the jacobian matrix at the NME point. The background color stands for the corresponding numerical $\omega$-limit of the principal trajectory~\eqref{case:interest:2:ic}, see Figure~\ref{fig:CoI3_highIRmut}, which are \textcolor{ColorNME}{NME}, \textcolor{ColorDFE}{DFE} and \textcolor{ColorCSE}{CSE}. The two vertical green lines indicate the region in which the system exhibits complex eigenvalues, \textit{i.e.}, where the jacobian has eigenvalues with non-zero imaginary parts, corresponding to the interval $\pi_1 \in (0.357.., 5.129..)$.}
    \label{fig:eigenvaluesNME}
\end{figure}

\noindent\textbf{$\bullet$ No master equilibrium (NME)} \\[1.1ex]
The spectrum of the jacobian matrix evaluated on the NME points
\begin{equation*}
\textrm{NME}:\quad(S^*,I_0^*,I_1^*,R^*,D^*)\times(g_0^*,g_1^*; v_0^*,v_1^*)=
\left(
\frac{\pi_1}{\beta_{11}}, 0 , 
\frac{1- \frac{\pi_1}{\beta_{11}}}{1+\frac{\pi_1}{\chi}},
\frac{\pi_1}{\chi} \cdot \frac{1- \frac{\pi_1}{\beta_{11}}}{1+\frac{\pi_1}{\chi}},
 0 \right)  \times   \left(0,1;0,\frac{\xi_1}{\gamma_1}\right),
\end{equation*}
defined for $\pi_1 < 40/7$, is denoted by
\begin{equation*}
\lambda_{\mathrm{NME}} = \Bigg\{-\frac{4}{5} ,\, -\frac{1}{2} ,\, -\frac{1}{2} ,\, -\frac{1}{5} ,\, -\frac{1}{5} ,\, 0 ,\, -\frac{54 - \sqrt{2}\psi}{7\,{\left(\mathrm{\pi_1}+2\right)}} ,\, -\frac{54 + \sqrt{2}\psi}{7\,{\left(\mathrm{\pi_1}+2\right)}}\,\Bigg\}
\label{eq:eigenvalues_NME}
\end{equation*}
with $\psi =\sqrt{\Delta}$ and $\Delta=\Delta(\pi_1):=49 \pi_1^3 -84 \pi_1^2 -924\pi_1+338$. As in the previous case, it has one zero eigenvalue (associated with the neutral (or central) manifold $I_0=0$), five real negative eigenvalues, and a couple of complex conjugate eigenvalues
\[
\lambda_{\mathrm{nme}^-} = -\frac{54 - \sqrt{2}\psi}{7\,{\left(\mathrm{\pi_1}+2\right)}}, \qquad \quad 
\lambda_{\mathrm{nme}^+} = -\frac{54 + \sqrt{2}\psi}{7\,{\left(\mathrm{\pi_1}+2\right)}},
\]
which can be real or complex depending on the sign of $\Delta$. This discriminant is $\Delta<0$ if $\pi_1\in (0.356.., 40/7)$ and $\Delta>0$ if $\pi_1<0.356..$ or $\pi_1>40/7$. Figure~\ref{fig:eigenvaluesNME} shows the real and imaginary part of these two eigenvalues as functions of $\pi_1$.

Therefore, for $\pi_1\in (\pi_0, 40/7)$, we have two regions of stability according to its value:
\begin{itemize}
\item[(a)] Case $\pi_1 \in (\pi_0, 5.128..)$: in this scenario, $\psi=\sqrt{\Delta}$ is complex and so $\lambda_{\mathrm{nme}^{\pm}}$ are a couple of complex conjugate eigenvalues, namely:
\[
\lambda_{\mathrm{nme}^{-}} = -\frac{54}{7(\pi_1+2)} + \rmi \frac{\sqrt{-2\Delta}}{7(\pi_1+2)}, \qquad \quad 
\lambda_{\mathrm{nme}^{+}} = \overline{\lambda_{\mathrm{nme}^{-}}}
\]
Since their real part is negative, the corresponding NME point is attractor (out of the neutral manifold). In the $2$-dimensional manifold associated to   
$\lambda_{\mathrm{nme}^{\pm}}$, the trajectories exhibit damped oscillations. It is known that, on this plane, and close to the equilibrium point, the quasi-period of such oscillations is
\begin{equation} 
T(\pi_1) \simeq \frac{2\pi}{\mathrm{Im}\, (\lambda_{\mathrm{nme}^{-}}) } = 
\frac{\sqrt{2} \pi (7\pi_1+2)}{\sqrt{-\Delta}},
\label{eq:quasiperiod}
\end{equation}

where $\Delta=\Delta(\pi_1)$ (Figure~\ref{fig:coi1:quasiperiod:NME}).
\begin{figure}[!ht]
\centering
\includegraphics[scale = 1]{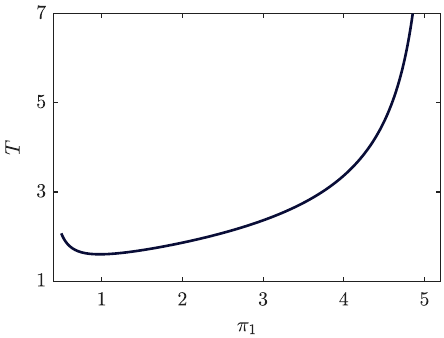}
\caption{Quasi-period $T(\pi_1)$, see~\eqref{eq:quasiperiod}, of the damped oscillations in a neighbourhood of the NME point for $\pi_1 \in (1/2, 5.128..)$.}
\label{fig:coi1:quasiperiod:NME}
\end{figure}

\item[(b)] Case $\pi_1 (5.128.., 40/7)$: now the discriminant $\Delta>0$ implies that both $\lambda_{\mathrm{nme}^{\pm}}$ are real. Furthermore, $\lambda_{\mathrm{nme}^{+}}<0$ and $\lambda_{\mathrm{nme}^{-}}<0$. Hence, the NME point is stable.
\end{itemize}

\begin{figure}[!ht]
    \centering
    \includegraphics[width = \textwidth]{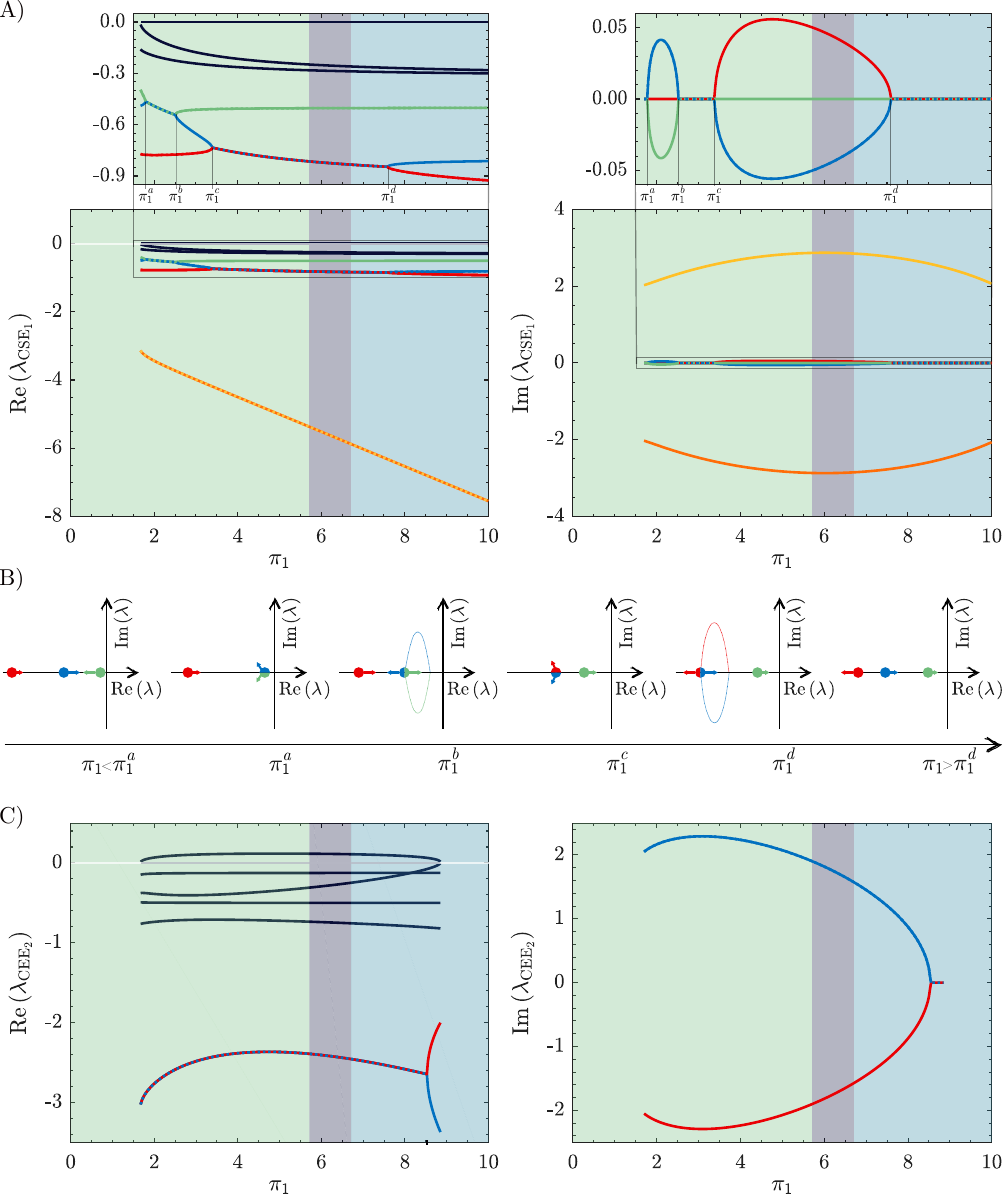}
    \caption{(A) Real (left) and imaginary (right) parts of the eigenvalues associated with the numerically computed equilibrium point $\mathrm{CEE}_1$. Insets enhance resolution within the grey-marked region. (B) Schematic representation of the evolution of the three eigenvalues undergoing bifurcations at $\pi_1 = \pi_1^j$, with $j = a,b,c,d$ (see main text). (C) Same as (A), for the equilibrium $\mathrm{CEE}_2$. In (A) and (C), background color indicates the $\omega$-limit of the principal trajectory (Figure~\ref{fig:CoI3_highIRmut}): \textcolor{ColorNME}{NME}, \textcolor{ColorDFE}{DFE}, or \textcolor{ColorCSE}{CEE}. Color coding links real and imaginary parts of each eigenvalue; black indicates a vanishing imaginary part.}
    \label{fig:eigenvaluesCSE}
\end{figure}

\noindent\textbf{$\bullet$ Co-circulating strains equilibrium (CSE)} \\[1.1ex]
For the set of parameters~\eqref{caseinterest1:CSE:parameters}
our model has two distinct families of CSE points, denoted as $\mathrm{CSE}_1$ and $\mathrm{CSE}_2$ and found numerically and illustrated in Figure~\ref{fig:CEEpoints}. To assess their local stability, we numerically compute their jacobian matrices and calculate their spectrum. Figure~\ref{fig:eigenvaluesCSE} shows the results of this analysis.

Figure~\ref{fig:eigenvaluesCSE}A shows the real and imaginary parts of the eigenvalues as functions of $\pi_1$ for the equilibrium point CSE$_1$ (left panel in Figure~\ref{fig:CEEpoints}), for a range of values of the parameter $\pi_1$. In this case, all eigenvalues exhibit negative real parts across this interval, meaning that the equilibrium is asymptotically stable. Nonetheless, the presence of complex conjugate eigenvalues introduces damped oscillatory modes into the local dynamics near the equilibrium. Two eigenvalues, highlighted in orange and yellow, remain a complex conjugate pair accross the entire range of $\pi_1$ explored. This persistent non-zero imaginary part indicates the existence of an oscillatory mode whose frequency and decay rate vary with $\pi_1$. In contrast, three other eigenvalues (blue, green, and red) alternate between being real and forming complex conjugate pairs. These transitions correspond to bifurcations occurring at critical values of $\pi_1$, denoted $\pi_1^a$, $\pi_1^b$, $\pi_1^c$, and $\pi_1^d$ in the figure. At these bifurcation points, the corresponding eigenvalues collide on the real axis and subsequently branch off into the complex plane (or \textit{vice versa}), thereby introducing or annihilating additional oscillatory modes. Figure~\ref{fig:eigenvaluesCSE}B provides a schematic representation of these bifurcations in the complex plane. The trajectories of the blue, green, and red eigenvalues are illustrated to emphasize their transitions between real and complex configurations. Collectively, these bifurcations enrich the system's dynamical landscape, enabling multiple qualitatively distinct oscillatory regimes depending on the parameter $\pi_1$.

Finally, Figure~\ref{fig:eigenvaluesCSE}C shows the eigenvalue spectrum corresponding to $\mathrm{CEE}_2$ (see the right panel in Figure~\ref{fig:CEEpoints}). Across the range of $\pi_1$ studied, one eigenvalue exhibits a strictly positive real part, which reveals the instability of $\mathrm{CEE}_2$. Additionally, two eigenvalues form complex conjugate pairs, indicating the presence of an oscillatory mode, although the dynamics would be ultimately dominated by the unstable direction. This complex conjugate pair also exhibits a bifurcation right before $\mathrm{CSE}_2$ vanishing. Due to the repelling nature of $\mathrm{CSE}_2$, the simulated principal trajectories will never converge to this equilibrium.

\subsubsection[The limit case pi1 to infty]{The limit case $\mathbf{\pi_1 \rightarrow +\infty}$}
\label{sec:coi1:limitcase}

\begin{figure}[!ht]
    \begin{minipage}{0.55\textwidth}
    \centering
    \includegraphics[scale = 1]{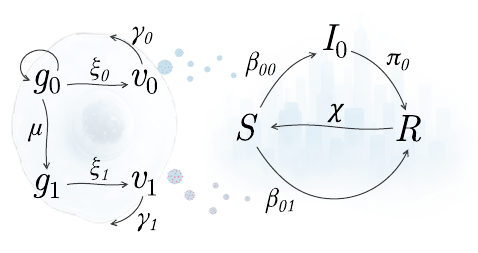}
    \end{minipage}
    \hfill
    \begin{minipage}{0.44\textwidth}
    \centering
    \hspace{-1.75cm}
    \begin{tabular}{c}
        Stoichiometry of the macroscopic model\\
        \hline\noalign{\vskip 4pt}
        $S+I_0\xrightarrow{\beta_{00}} I_0+I_0$\\
        $S+I_0\xrightarrow{\beta_{01}} I_0+R$\\
        $I_0\xrightarrow{\pi_0}R$ \\
        $R \xrightarrow{\chi} S$\\[3pt]
        \hline
    \end{tabular}
    \end{minipage}
    \caption{(left) Schematic representation of the reduced model in the limit case $\pi_1\rightarrow\infty$. (Right) Macroscopic reaction netowrk of the model. Note in particular the reaction in which contact between $I_0$ and $S$, mediated by mutant virions $v_1$ (with transmission rate $\beta_{01}$), results in an individual that is neither infected nor susceptible, \textit{i.e.}, belonging to the recovered class $R$. The microscopic component remains unchanged, except that $\tilde{f}_1 = 0$ as $I_1 \to 0$, implying that no substrate is available for mutant replication.}
    \label{fig:reducedModel}
\end{figure}

An extreme scenario arises when the mutant strain induces a very rapid recovery in the host, $\pi_1\gg 1$. This regime can be mathematically approximated by taking the limit $\pi_1 \rightarrow +\infty$. In this extreme case, the dynamics resemble a situation in which the master variant gives rise to a mutant strain that lacks self-replicative capacity (see Figure~\ref{fig:reducedModel}). In this setting, the mutant can only increase in frequency due to erroneous replication of the master genome ($g_0$), which occasionally produces virions ($v_1$) capable of infecting new hosts. However, upon infection by this defective strain, the host rapidly recovers and temporarily acquires immunity, governed by the parameter $\chi$.

Under these conditions, we may assume that $I_1 \rightarrow 0$, since recovery is effectively instantaneous. This implies $\nu_1 \rightarrow 0$ and either $\nu_0 \rightarrow 1$ (when $I_0 \ne 0$) or $\nu_0 \rightarrow 0$ (if $I_0=0$). Consequently, $\beta_{11} \rightarrow 0$, $\tilde{f}_1=f_1 \nu_1 \rightarrow 0$ and either $\tilde{f}_0=f_0 \nu_0 \rightarrow f_0$ (if $I_0\ne 0$) or $\tilde{f}_0=f_0 \nu_0 \rightarrow 0$ (if $I_0=0$). Henceforth, we assume that all variables and parameters take their corresponding limiting values.

From Section~\ref{se:caseofinterest1:equilibria}, we know that two type of equilibrium points exist in this limit: the DFE and the CSE. In the former case, $I_0=0$ implies $\nu_0=0$ and $\tilde{f}_0=0$, leading to a trivial genomic system $\dot{g}_0=0$, $\dot{g}_1=0$. Its equilibrium solutions are of the form $(g_0^*, 1-g_0^*)$ for $g_0^*\in [0,1]$. Together with the corresponding virion equilibria, these states yield to the previously defined $\QSmicro$ (see Proposition~\ref{prop:DFE}). In this case, it is straightforward to show that the unique macroscopic DFE is $(S^*, I_0^*, R^*) = (1, 0, 0)$. In contrast, in the CSE case, when $I_0 > 0$, the microscopic system becomes
\[
\dot{g}_0 = f_0g_0(1-\mu-g_0), \qquad \dot{g}_1 = f_0g_0 (\mu - g_0),
\]
for $0<\mu<1$. This system has two equilibrium points: $(g_0^*,g_1^*)=(0,1)$ (repeller) and $(g_0^*, g_1^*)=(1-\mu, \mu)$ (attractor), which derive into
\[
(v_0^*,v_1^*)=\left( 0, \frac{\xi_1}{\gamma_1} \right) \qquad \textrm{and} \qquad (v_0^*,v_1^*) = \left( \frac{\xi_0}{\gamma_0}(1-\mu), \frac{\xi_1}{\gamma_1} \mu \right),
\]
respectively, for the virion equilibria. In this scenario, the macroscopic system reads
\[
\dot{S}=- \left( \beta_{00}+\beta_{01} \right) I_0 S + \chi R, \qquad 
\dot{I}_0 = \left( \beta_{00} S- \pi_0 \right) I_0, \qquad 
\dot{R} = \beta_{01} I_0 S + \pi_0 I_0 - \chi R.
\]
Notice that $v_0^*=0$ implies $\beta_{00}=0$ and so, from the second equation above, that $\pi_0=0$ (since $I_0>0$), a contradiction with the fact that we are assuming $\pi_0=2$ in this section. 
The second microscopic equilibrium solution is $\QSmicromu$,
\[
(g_0^*, g_1^*; v_0^*, v_1^*) = \left( 1-\mu, \mu; \frac{\xi_0}{\gamma_0}(1-\mu), \frac{\xi_1}{\gamma_1} \mu \right)  
\]
already defined in~\eqref{b:mu:equal:1}. The corresponding equilibrium for the macroscopic system is easily derived:
\begin{equation}
S^* = \frac{\pi_0}{\beta_{00}}, \qquad I_0^* = \frac{\chi}{\chi + \left( 1 + \dfrac{\beta_{01}}{\beta_{00}}\right) \pi_0}, \qquad 
R^* = 1 - S^* - I_0^*,
\label{CSEeq:points:limit:case}
\end{equation}
such that the three variables sum $1$.

Concerning their local stability, the computation of the spectra of their associated jacobian matrices, gives rise to:
(a) $\{ -\chi_1, \beta_{00} - \pi_0 \}$ for the DFE point $(1,0,0)$, which means that it is stable if $\pi_0 > \beta_{00}$ (\textit{i.e.}, it has $\mathcal{R}_0<1$) and unstable otherwise;
(b) for the CSE point~\eqref{CSEeq:points:limit:case}, it is straightforward to show that its determinant is positive, the trace negative and, hence, it is always stable. Furthermore, its behavior moves between a stable node or focus (so, oscillations appear) depending on the values of $\beta_{00}$, $\beta_{01}$, $\pi_0$, and $\chi$. These parameters also determine, in the cases of bistability, the configuration of the corresponding basins of attraction.
\renewcommand{\arraystretch}{1.1}

Biologically, this setting is reminiscent of defective viral genomes (DVGs) which, despite their inability to autonomously sustain a replication cycle by themselves, can still be packaged into master-encoded viral particles and transmitted \cite{Vignuzzi2019}. Although such transmission does not establish a productive infection, it may nevertheless be sufficient to trigger an immune response in the host, for example through superinfection exclusion or the induction of immune memory. Indeed, this property has opened the tantalizing possibility of using engineered versions of DVGs as self-transmissible vaccines \cite{Notton2014}.


\subsection{Case of interest 2: the burnout viral strain.}
\label{sec:burnoutVirus}

Below, we highlight three well‑documented examples that lie near the burnout corner of our parameter space, characterized by very high virulence ($\delta_1 \gg 0$), transmission rate sufficient to ignite large outbreaks (effectively large $a_1$), and slow recovery at the population level (effectively small $\pi_1$). Each example exhibits the core macroscopic level signature reproduced by our model: rapid epidemic growth, high mortality, and self‑limitation throughout host depletion or behavioral change, sometimes followed by longer‑term evolutionary relaxation.

Firstly, the introduction of myxoma virus (MYXV) into Australia to control invasive European rabbits provides a classic illustration~\cite{Kerr2012}. Initial mortality reached approximately 99 \%, and MYXV spread explosively via arthropod vectors, collapsing rabbit populations from an estimated $\sim600$ million to $\sim100$ million within a few years, a textbook burnout pulse (huge realized $a_1$ in a highly connected host–vector system combined with extreme $\delta_1$). Remarkably, subsequent coevolution reduced virulence and increased host resistance, gradually shifting the system away from the burnout regime~\cite{Kerr2013}. Secondly, the 2014-2016 West African outbreak of Zaire Ebolavirus (ZEBOV) was the largest filovirus outbreak on record, with estimated case fatality rates around 60 – 70\%~\cite{Kucharski2014}. Simple growth models during early spread gave $R_0 \approx 1.6 - 2.0$, indicating sustained human‑to‑human transmission \cite{Althaus2015}. These outbreaks typically burn out once susceptible depletion, improved infection prevention, and control reduce effective transmission, consistent with a high $\delta_1$ process that self‑limits in the absence of continuous replenishment of susceptibles.  Thirdly, since 2021 – 2023, highly pathogenic avian influenza virus (HPAIV) H5N1 clade 2.3.4.4b has caused unprecedented die‑offs in wild birds and mass mortality in marine mammals~\cite{Leguia2023}. These events reveal intense transmission in dense breeding colonies coupled with extremely high mortality, exemplifying the “ignite fast, die fast” pattern. Genomic epidemiology provides evidence consistent with mammal‑to‑mammal transmission in pinnipeds. In such colonies, the realized $a_1$ becomes very large while $\delta_1$ is extreme, leading to rapid exhaustion of susceptible and epizootic burnout~\cite{GamarraToledo2023}.

To explore the evolutionary consequences of virulence emergence, we devote this section to analyzing a scenario in which an initially avirulent master strain ($\delta_0 = 0$) with moderate transmission ($a_0 = 2$) gives rise to a mutant variant with higher virulence and transmissibility, here $\delta_1 = 1$ and $a_1 = 3$. All the remaining parameters of the model have been conveniently fixed at the following values:
\begin{eqnarray}
&&\chi = 2,\quad \delta_0 =0, \quad \delta_1=1, \quad \pi_0=\pi_1=0.2,
\quad  a_0=2, \quad a_1=3, \quad b_0=b_1=0.5, \nonumber \\
&&f_0 = 1,\quad  f_1=0.1, \quad \gamma_0=\gamma_1=0.5, \quad \xi_0 =\xi_1=3,\quad  \gamma_0 = 0.8. \label{caseinterest2:parameters} 
\end{eqnarray}
Again, the slow-fast time ratio is taken as $\eps = 0.01.$

\begin{figure}[!ht]
    \centering
    \includegraphics[width = \textwidth]{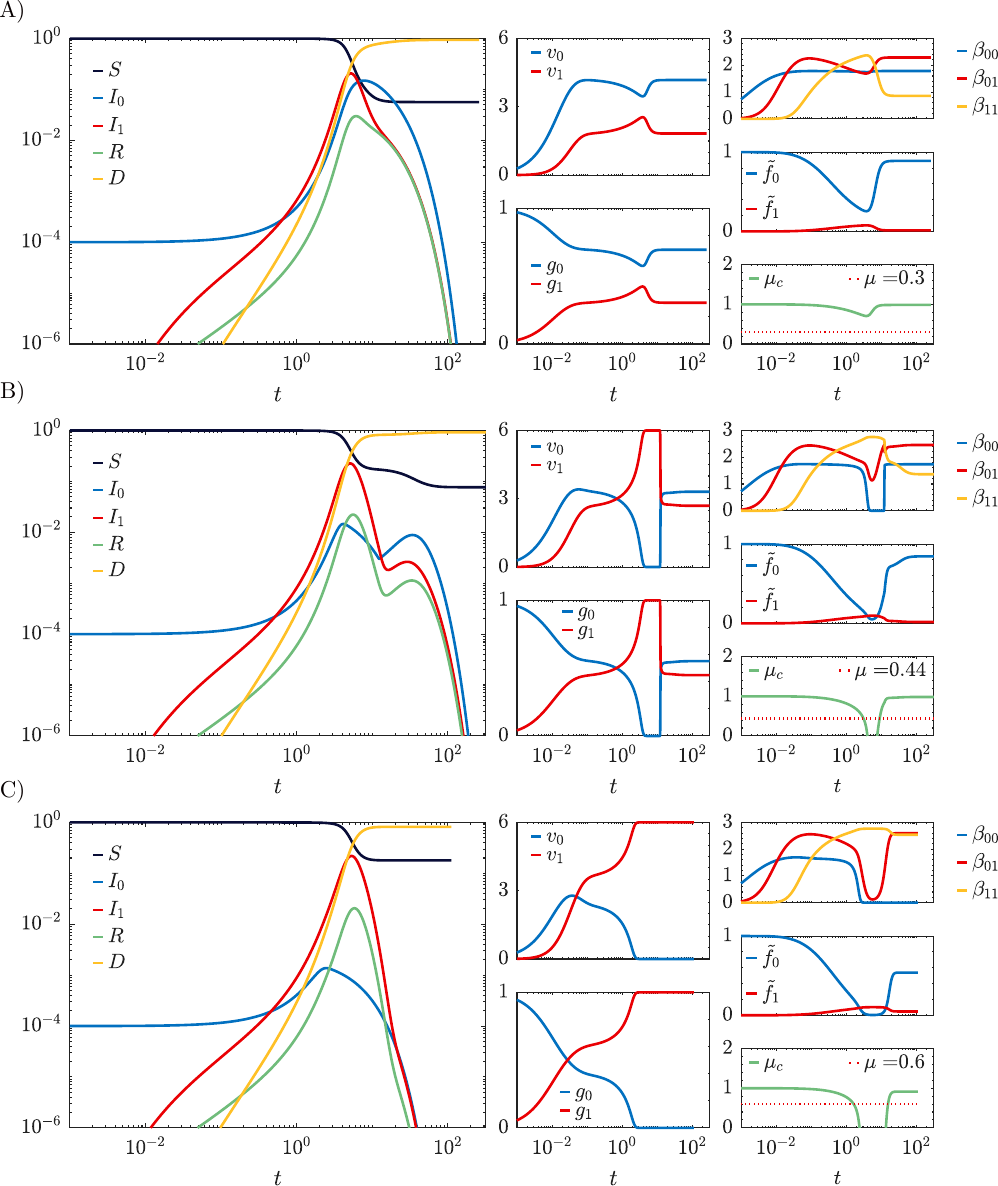}
    \caption{Numerical integrations of the system for non-varying parameter set \eqref{caseinterest2:parameters}. The slow to fast time ratio is fixed at $\eps = 0.01$. To avoid false recoveries from pseudo-extinction, we set to zero any variable below the numerical tolerance of the ODE numerical integrator (\emph{i.e.} $10^{-14}$).}
    \label{fig:CoI2_HorrorSupermarket}
\end{figure}

As a first step in this study, we vary the mutation probability $\mu$ and analyze its impact on the dynamics of the system. In particular, we monitor: (\textit{i}) the evolution of the state variables; (\textit{ii}) the effective transmission rates of both strains ($\beta_{00}, \beta_{01}$, and $\beta_{11}$); and (\textit{iii}) the instantaneous fitness of each variant ($\tilde{f}_0$ and $\tilde{f_1}$) and the dynamic critical mutation threshold ($\mu_c=\mu_c(t)$), both characteristics being crucial in the inter-level connection. Three snapshots of this analysis are shown in Figure~\ref{fig:CoI2_HorrorSupermarket}.

We highlight the particularly interesting case arising in the central panels of Figure~\ref{fig:CoI2_HorrorSupermarket}, where the master variant appears to overcome a pseudo-error catastrophe. This abrupt transition arises due to the micro– to macroscopic coupling inherent in the model.

When the mutation probability is very high (\textit{e.g.}, $\mu \gtrsim 0.4$), continuous generation of mutant genomes $g_1$ drives a rapid rise in virions $v_1$, which boosts the strain‑specific transmission rate and shifts prevalence so that $I_1 \gg I_0$. This macroscopic level shift depresses the master’s context‑dependent error threshold $\mu_c(t)$, pushing the within‑host quasispecies beyond the threshold and causing a transient microscopic pseudo‑extinction of $g_0$, a behavior consistent with error‑catastrophe experiments in poliovirus exposed to ribavirin~\cite{Crotty2001}. Because the mutant is highly virulent $\delta_1 \gg 0$, mortality rapidly depletes $I_1$, transmission collapses, and conditions that favor the master return: as $g_1$ declines, the master’s effective replicative fitness $\tilde{f}_0$ increases, $\mu_c(t)$ rises above $\mu$, and residual $g_0$ genomes can resume replication. A macroscopic scale analogue of this “ignite fast, die fast” endgame appears in the aforementioned HPAIV wildlife pulses, where explosive spread with extreme lethality in dense colonies (\emph{e.g.}, sea lion mass mortality in Peru) is followed by rapid chain collapse \cite{GamarraToledo2023}. Ultimately, the epidemic wanes as cumulative deaths disrupt infection chains; in our parameter region with $\delta_0=0$ and $\delta_1>0$, only disease‑free equilibria (DFE) are feasible, and numerical simulations show the principal trajectory converges to DFE across comparable parameter sets. We analyze the DFE behavior as a function of $\mu$ for the representative parameter choice given in \eqref{caseinterest2:parameters}.

\bigskip

\begin{figure}[!ht]
\centering
\includegraphics[scale = 0.85]{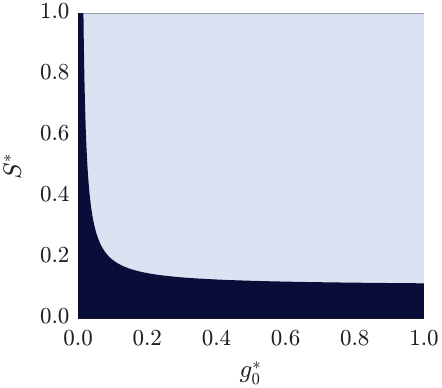}
\caption{Stability regions of the DFE plane. The attractive (unstable) region of each DFE point is represented in black (grey) colour, in terms of $g_0^*$ and $S^*$. The hyperbola $\psi(S^*,g_0^*)=0$ dividing both zones is given by~\eqref{caseofinterest2:DFE:vap:h}.} 
\label{fig:COI2_DFEstability}
\end{figure}

\noindent \textbf{$\bullet$ Disease-free equilibria (DFE).}
From Proposition~\ref{prop:DFE}, they are of the form:
\begin{equation*}
    (S, I_0, I_1, R, D) \times (g_0, g_1 \, ; v_0, v_1)=
    \left(S^*, 0,0,0, 1-S^*\right) \times \left(g_0^*, 1-g_0^* \,; \frac{g_0^*\xi_0}{\gamma_0},\frac{\left(1-g_0^* \right)\xi_1}{\gamma_1}\right).
\label{eq:HSTG:DFE}
\end{equation*}
These DFE equilibria fill a $2$-dimensional plane, governed by the variables $(S^*,g_0)$. Their local stability can be approached by their linearised system. Indeed, the spectrum of the corresponding 
jacobian matrices is given by
\begin{equation*}
\lambda_{\mathrm{DFE}} = \left\{ 0,\,0,\,0,\,0,\, -\chi ,\, -\gamma_0 ,\, -\gamma_1 ,\, -(\delta_1+\pi_1) ,\, \psi(S^*,g_0^*) \right\}.
\label{eq:eigenvalues_DFE}
\end{equation*}
where
\[
\psi(S^*,g_0^*)=\frac{a_0\xi_0g_0^* S^*}{b_0 \gamma_0 + \xi_0 g_0^*} - \pi_0.
\]
The four zero eigenvalues come from the first integrals $I_0=0$, $I_1=0$, $R=0$, and $S+I_0+I_1+R+D=1$. There are four more negative eigenvalues and a last one, $\psi(S^*,g_0^*)$, whose sign determines their (transversal) stability. Indeed, $\psi(S^*,g_0^*)=0$ is equivalent to the hyperbola (see Figure~ \ref{fig:COI2_DFEstability}).
\begin{equation}
S^* = \frac{\pi_0}{a_0} \left( 1 + \frac{b_0 \gamma_0}{\xi_0} \, \frac{1}{g_0^*}  \right), \qquad \quad g_0^* \in (0,1].
\label{caseofinterest2:DFE:vap:h}
\end{equation}


\section{Discussion}

We have presented a minimal, mechanistically grounded multiscale framework that links within-host quasispecies dynamics to population-level SIRS epidemiology through explicit, bidirectional coupling: transmission rates $\beta_{ij}$ depend on the instantaneous virion abundance $v_j(t)$ (microscopic to macroscopic), and the effective replicative rates $\tilde{f}_j$ are weighted by the prevalence of each infected host class (macroscopic to microscopic). The imposed slow-fast structure ($\varepsilon \ll 1$) yields a natural quasi-steady reduction of the genome-virion subsystem and explains why, on epidemiological timescales, transmission can be treated as time-varying process inherited from the fast layer while preserving a conserved quasispecies backbone. This construction moves beyond \textit{ad hoc} scale-bridging approaches by keeping the error-threshold mechanism explicit and intrinsically context-dependent.

Our NGM calculation at DFE yields $\mathcal{R}_0 = \max \left\{\beta_{00}/(\pi_0 + \delta_0), \beta_{11}/(\pi_1 + \delta_1)\right\}$, clarifying that the cross-seeding pathway $\beta_{01}$ does not affect invasion, because the jacobian block is triangular at DFE (no mutant virions are present). In contrast, $\beta_{01}$ plays a crucial role after invasion: once any $I_0>0$ exists, the $\beta_{01}$ channel accelerates mutant establishment and shapes transient dynamics of strain replacement. The slow-fast closure further justifies the introduction of time-dependent effective reproduction numbers $\mathcal{R}^{(0)}_t$ and $\mathcal{R}^{(1)}_t$, and the growth diagnostic $\mathcal{G}_t>0$ (see \eqref{def:G_t}) that decomposes overall expansion into outflow-weighted excess reproduction of each strain (Section~\ref{se:growing_pandemic}). These quantities admit clear biological interpretation and naturally connect to empirical data through longitudinal proxies of viral load $v_j(t)$, which control the saturating transmission functions $\beta_{ij}$.

A central theoretical contribution of this work is the demonstration that the error threshold becomes epidemic-state dependent (see \eqref{def:genom:critical:mutation:rate}). As a transmissible mutant rises in prevalence ($\nu$ decreases), the critical mutation rare $\mu_c(\nu)$ decreases, so that, at fixed polymerase fidelity, the master $g_0$ can transiently fall below its error threshold. This produces a pseudo-error catastrophe driven by cross-scale feedback rather than a biochemical change in replication fidelity (Section \ref{sec:burnoutVirus}; central panel of Figure~\ref{fig:CoI2_HorrorSupermarket}). As mortality and host depletion subsequently collapse $I_1$, $\mu_c(\nu)$ increases again above $\mu$, allowing residual $g_0$ to resume replication, thereby explaining the observed ``disappear-then-reappear'' microscopic signature during a macroscopic burnout pulse. This mechanism cannot be recovered when within- and between-host scales are analyzed in isolation.

The model yields sharp feasibility inequalities that delineate endpoints and coexistence regimes. For example, at candidate coexistence one obtains the condition $\pi_0\,\beta^{*}_{00} < \pi_1\,\beta^{*}_{11}$ (Section~\ref{sec:eqPoints}), where starred quantities denote the quasi-steady transmission rates inherited from the fast layer. Because $\beta^{*}_{ij}$ retain hyperbolic saturation and depend explicitly on $(\xi_j,\gamma_j,f_j)$, these inequalities make transparent how within-host production efficiencies bias population-level toward coexistence, no master or disease-free endpoints. 

Linear stability analysis reveals complex conjugate eigenvalues for NME and CSE across broad parameter regions (Figures~\ref{fig:eigenvaluesNME} and~\ref{fig:eigenvaluesCSE}), predicting damped oscillations in $I_0,I_1$ and $g_j,v_j$, even in the absence of external forcing, a distinctive dynamical signature of the multiscale feedback. The quasi-period near NME follows directly from this complex eigenvalue pair (see~\eqref{eq:quasiperiod} and Figure~\ref{fig:coi1:quasiperiod:NME}).

Two illustrative regimes bracket the model's behavior and directly connect to well-documented biology. In the ``vaccine-like'' quadrant ($\delta_0=\delta_1=0$, $\pi_1\gg \pi_0$), the mutant benignly replaces and persists at appreciable prevalence, generating high seroprevalence without excess mortality (Section~\ref{sec:COI1}; Figure~\ref{fig:CoI3_highIRmut}), consistent with live attenuated, strongly immunizing circulation. In the ``burnout'' quadrant ($\delta_1\gg0$, modest or weak $\pi_1$), a hypervirulent mutant self-limits throughout host depletion: the epidemic surge can pass through a pseudo-extinction of $g_0$ before collapsing back to the DFE, as observed numerically across a wide range of mutation rates (Section~\ref{sec:burnoutVirus} and Figure~\ref{fig:CoI2_HorrorSupermarket}). Together, these extremes emphasize that epidemic composition, rather than intrinsic polymerase fidelity alone, governs transient excursions across error thresholds.

We  can also draw several practical implications from our study. (\textit{i}) Transient sequencing signatures: during mutant ascents, we predict a temporary loss of the master sequence in within-host spectra, coupled to peaks in incidence, interpretable as a population-driven passage through the threshold $\mu_c(\nu)$. (\textit{ii}) Oscillatory returns: damped endemic waves should arise where the two-strain slow system has complex eigenvalues; these oscillations can be quantitatively fitted with virion-dependent transmission functions $\beta_{ij}$ applied to joint incidence and viral load time series. (\textit{iii}) Design criteria for benign replacement: the coexistence inequalities translate empirically measured load-transmission curves $(a_j,b_j)$ and immune induction $(\pi_1)$ into operational thresholds for safe dominance.

We intentionally compressed within-host diversity into a two-type, unidirectional quasispecies and neglect explicit immune kinetics, host heterogeneity, demographic turnover (births and immigration of susceptibles), contact structure, and stochastic fade-out. These simplifying choices enable analytical tractability but preclude phenomena such as diversification fronts, antigenic escape, and age or network effects on epidemic thresholds. The explicit slow-fast structure of the model suggests a rigorous Tikhonov reduction to a closed slow system, a global bifurcation analysis (saddle-node, Hopf) expressed in terms of measurable parameters (\emph{e.g.}, peak viral load and decay), and extensions incorporating reverse mutation, multi-lineage competition, host heterogeneity in susceptibility and waning immunity, and stochastic invasion dynamics. Finally, the $\pi_1 \rightarrow \infty$ limit (Section~\ref{sec:coi1:limitcase}) motivates the study of defective viral genomes within the same multiscale closure, in which the cross-seeding pathway $\beta_{01}$ functions as an “immunizing contact” that induces host immunity without productive mutant infection.


\section{Concluding remarks}

By tying transmission directly to intracellular production and allowing epidemiological prevalence to feedback on effective replication, our model shows how selection pressures at one scale can induce phase transition-like shifts at the other. Most notably, how a mutant surge can transiently trigger an intrahost error catastrophe and how extreme virulence self-limits transmission and spread. The ``vaccine-like'' and ``burnout'' regimes bracket a continuum of realistic outcomes, unify microcroscopic and macroscopic thresholds (from $\mu_c(\nu)$ to $\mathcal{R}_0$), and yield testable signatures in coupled incidence–viral-load time series. We anticipate that combining quasi-steady reductions, joint inference on $\beta_{ij}(v_j)$ from viral load and contact data, and explicit heterogeneity will transform this conceptual bridge into a practical tool for anticipating which mutants are poised to dominate, coexist, or extinguish themselves—and why.

\section{Acknowledgements}

JCM-S was supported by Generalitat Valenciana grant ACIF/2021/296. JTL and JS thank the AEI, through the Mar\'ia de Maeztu Program for Units of Excellence in R\&D (CEX2020-001084-M) and the Generalitat de Catalunya CERCA Program for institutional support. JTL has been also funded by grants PID2021-122954NB-I00 and PID2024-155942NB-I00 funded by MICIU/AEI/10.13039/501100011033 and ``ERDF a way of making Europe". SFE was supported by grants PID2022-136912NB-I00 funded by MCIU/AEI/10.13039/501100011033 and by ``ERDF a way of making Europe”, and CIPROM/2022/59 funded by Generalitat Valenciana. JCM-S, JTL and SFE thank support from the Santa Fe Institute.

\bibliographystyle{unsrt}
\bibliography{bibliography}

\newpage
\begin{center}
{\Large{\textsc{Appendix}}    }
\end{center}

\appendix
\setcounter{equation}{0}
\renewcommand{\theequation}{\thesection\arabic{equation}}
\setcounter{figure}{0}
\renewcommand{\thefigure}{\thesection\arabic{figure}}
\setcounter{table}{0}
\renewcommand{\thetable}{\thesection\arabic{table}}

\setcounter{figure}{0}

\section{Equilibrium points: proofs of the propositions}
\label{sec:proofs:eqPoints}

\subsection{Proof of Lemma~\ref{lemma:micro}}

Notice that if one defines $\tilde{f}_0=f_0 \nu^*$ and $\tilde{f}_1=f_1 \nu^*$, the microscopic system becomes
\[
\eps \dot{g}_0 = \tilde{f}_0 (1-\mu) g_0 - \Phi g_0, \qquad
\eps \dot{g}_1 = \tilde{f}_0 \mu g_0 + \tilde{f}_1 g_1 - \Phi g_1,
\]
with $\Phi=\tilde{f}_0 g_0 + \tilde{f}_1 g_1$. This is an standard quasispecies model with fitnesses $\tilde{f}_0$ and $\tilde{f}_1$. It is well-known \cite{Eigen1977} that there are two possible scenarios for its equilibrium points in terms of the relation between the mutation probability $\mu$ and the corresponding critical mutation probability.
\[
\mu_c^* = 1 - \frac{\tilde{f}_1}{\tilde{f}_0} = 1 - \frac{f_1}{f_0} \left( \frac{1}{\nu^*} - 1 \right).
\]
\begin{itemize}
 \item[(i)] If $0<\mu < \mu_c^*$ then the equilibrium point is
\[
(g_0^*,g_1^*) = \left( 1 - \frac{\mu}{\mucrit^*}, \frac{\mu}{\mucrit^*} \right).
\]
\item[(ii)] On the contrary, if $\mu_c^* \leq \mu \leq 1$, the equilibrium point is uniquely composed by mutant genomes: $(g_0^*,g_1^*) =  (0,1)$.
\end{itemize}
The expressions for the equilibrium virions derive straightforwardly from the ones for the genomes.
\qed

\subsection{Proof of Proposition~\ref{prop:DFE} (DFE)}
Under the assumption $I_0^*=I_1^*=0$ the macroscopic system becomes
\[
\dot{S}=\chi R, \quad \dot{I}_0=0, \quad \dot{I}_1=0, \quad \dot{R}=-\chi R, \quad
\dot{D}=0.
\]
Equating them to zero (to seek for equilibrium points) we get $R^*=0$ if $\chi\ne 0$ or arbitrary $R$ when $\chi=0$. Regarding the microscopic system, since $I_0^*=I_1^*=0\Rightarrow \nu_0^*=\nu_1^*=0$ we have $\tilde{f}_0=\tilde{f}_1=\Phi=0$ and so
\[
\eps \dot{g}_0 =0, \qquad
\eps \dot{g}_1 =0,
\]
which leads to equilibria with arbitrary $(g_0, g_1)$ such that $g_0+g_1=1$. Regarding virions, it is straigthforward to show that the equilibria are of type
\[
v_0=\frac{\xi_0}{\gamma_0} g_0, \qquad v_1=\frac{\xi_1}{\gamma_1} g_1, 
\]
for any $(g_0, g_1)$ satisfying $g_0+g_1=1$.

\qed

\subsection{Proof of Proposition~\ref{prop:NME} (NME)}

Substituting $I_0^*=0$ in $\dot{D}=\delta_0 I_0^* + \delta_1 I_1^*=0$ it turns out that $\delta_1 I_1^* =0$. Since $I_1^*>0$ it follows that, necessarily, $\delta_1=0$. We consider two cases:

\noindent (\textit{i}) Case $\beta_{11}>0$. This is divided, in its turn, into:
\begin{itemize}
\item[(\textit{i$_1$})] Case $\chi>0$. From $\nu^*=0$ we have, in particular, that $\beta_{01}=0$.
Substituting $I_0^*=0$ into~\eqref{eq:sirs:I1}, $\dot{I}_1=0$,  we obtain
\[
\left( \beta_{11} S - (\pi_1 + \delta_1) \right) I_1 =0
\Rightarrow
S= \frac{\pi_1 + \delta_1}{\beta_{11}},
\]
which becomes
\begin{equation*}
S^*= \frac{\pi_1}{\beta_{11}}.
\label{eqpoints:S:2}
\end{equation*}
since $\delta_1=0$. Substituting now $I_0=0$ into~\eqref{eq:sirs:R} we get
\begin{equation*}
\dot{R}=\pi_1 I_1 - \chi R =0 \Rightarrow R^*=\frac{\pi_1}{\chi} I_1,
\label{eqpoints:2:I1}
\end{equation*}
which is well defined since $\chi>0$. The variable $D^*$ is defined such that  $S^*+I_1+R^*+D^*=1$.

\item[(\textit{i$_2$})]  Case $\chi=0$: we have, from the same argument as above, that $S^*=\pi_1/\beta_{11}$. However,
\[
\dot{R}=0 \Leftrightarrow \pi_1 I_1^* = 0 \Rightarrow \pi_1 = 0 \Rightarrow S^*=0.
\]
That is, in this scenario, necessarily $\pi_1=0$ and the equilibrium point for the macroscopic system is of the form $(S^*,I_0^*, I_1^*, R, D)=(0,0,I_1,R,D)$ with arbitrary $R$, $D$, such that $I_1^* + R + D =1$.
\end{itemize}

\medskip

\noindent (\textit{ii}) Case $\beta_{11}=0$: the conditions for having an equilibrium point reduce to $\chi R=0$ and $ -\pi_1 I_1^* =0$. From the second equation it follows that, necessarily, $\pi_1=0$.  From the first one, we have two possible solutions:
\begin{itemize}
\item[(\textit{ii$_1$})] If $\chi>0$: then $R^*=0$ and the point are of the form $(S,I_0^*, I_1^*, R^*, D)=(S,0,I_1,0,D)$ with arbitrary $S$ and $D$ satisfying that $S+I_1^*+D=1$.
\item[(\textit{ii$_2$})] If $\chi=0$:  the equilibrium is of the form $(S,I_0^*, I_1^*, R, D)=(S,0,I_1,R,D) $ with arbitrary $R$, $S$ and $D$ such that $S+I_1^*+R+D=1$.

\end{itemize}

\medskip

In both two cases (\textit{i}) and (\textit{ii}),  at the genome's level, $\nu^*=0$ translates into $\Phi=f_1 g_1$ and so
\[
\eps \dot{g}_0 = - \Phi g_0 = - f_1 g_1 g_0 =0, \qquad
\eps \dot{g}_1 = f_1 g_1 - \Phi g_1 = f_1 g_1(1-g_1)=0.
\]
Since $g_0+g_1=1$ the only two solutions of the latter equations is either $(g_0^*,g_1^*)=(1,0)$ or $(g_0^*,g_1^*)=(0,1)$.  Therefore, joint to the corresponding solutions for the virion's system, one obtains the following two possible equilibrium points for the microscopic system:
\[
\QSmicro^{(0)}: \ (g_0^*,g_1^*; v_0^*, v_1^*) = \left( 1,0, \frac{\xi_0}{\gamma_0}, 0 \right)
\qquad \textrm{and} \qquad
\QSmicro^{(1)}: \ (g_0^*,g_1^*; v_0^*, v_1^*) =\left( 0,1, 0, \frac{\xi_1}{\gamma_1} \right).
\]
For $\QSmicro^{(1)}$ the condition $\beta_{11}=0$ implies that either $a_1=0$ or $\xi_1=0$. Condition $\beta_{11}>0$ is not compatible with equilibrium $\QSmicro^{(0)}$ since $\beta_{11}>0 \Rightarrow v_1^*>0$.
\qed

\subsection{Proof of Proposition~\ref{prop:NmutE} (NmutE)}
Let us start computing the possible equilibrium points of the microscopic system. Indeed, $I_1^*=0 \Rightarrow \nu^*=1 \Rightarrow \Phi=f_0g_0$. Thus, the equilibrium points  of the genome's system must solve
\begin{eqnarray*}
\eps \dot{g}_0 &=& f_0(1-\mu) g_0 - (f_0g_0)g_0 = f_0 g_0 (1-\mu - g_0)=0 \\
\eps \dot{g}_1 &=& f_0 \mu g_0 - (f_0g_0) g_1 = f_0g_0(\mu - g_1) =0.
\end{eqnarray*}
Since $f_0>0$ they become
\[
g_0(1-\mu -g_0)=0 \qquad \textrm{and} \qquad g_0(\mu - g_1) =0,
\]
which leads to two possible cases depending on $g_0$ vanishing or not. Thus,
\begin{itemize}
 \item[(\textit{a})]  Case $g_0=0$: so $g_1=1$ and then $(g_0^*, g_1^*)=(0,1)$. Regarding the virion's system, it becomes
 \begin{eqnarray*}
 \eps \dot{v}_0 &=& \xi_0 g_0^* - \gamma_0 v_0 = - \gamma_0 v_0 = 0 \Rightarrow v_0^*=0. \\
 \eps \dot{v}_1 &=& \xi_1 g_1^* - \gamma_1 v_1 = \xi_1 - \gamma_1 v_1 =0 \Rightarrow v_1^* = \frac{\xi_1}{\gamma_1}.
 \end{eqnarray*}
So, the equilibrium point is
\[
 (g_0^*, g_1^*; v_0^*, v_1^*) = \left( 0,1,0,\frac{\xi_1}{\gamma_1} \right)=:\QSmicro^{(1)}.
\]
Notice that $v_0^*=0$ implies that $\beta_{00}^*=0$ and from $\nu^*=1$ it follows that
$\beta_{11}^*=0$.

\item[(\textit{b})]
Case $g_0\ne 0$: this implies $(g_0^*, g_1^*) = (1-\mu, \mu)$ with $0<\mu \leq  1$. Again, substituting into the virion's  system we get
\begin{eqnarray*}
 \eps \dot{v}_0 &=& \xi_0 g_0^* - \gamma_0 v_0 = 0 \Rightarrow v_0^*= \frac{\xi_0}{\gamma_0}(1-\mu), \\
 \eps \dot{v}_1 &=& \xi_1 g_1^* - \gamma_1 v_1 =0 \Rightarrow v_1^* = \frac{\xi_1}{\gamma_1} \mu.
\end{eqnarray*}
So, the equilibrium point for the microscopic system is
\begin{equation}
(g_0^*, g_1^*; v_0^*, v_1^*) = \left( 1-\mu,\mu, \frac{\xi_0}{\gamma_0}(1-\mu), \frac{\xi_1}{\gamma_1} \mu \right) =: \QSmicromu.
\label{b:mu:equal:1}
 \end{equation}
 Moreover, $\beta_{11}=0$ (since $\nu^*=1$) and
 \[
 \beta_{00}^* = \frac{a_0 v_0^*}{b_0 + v_0^*}, \qquad
 \beta_{01}^* = \frac{a_1 v_1^*}{b_1+v_1^*},
 \]
which vanish if and only if $a_0=0$ or $v^*_0=0$ and $a_1=0$ or $v_1^*=0$, respectively.
 \end{itemize}

Having in mind these two equilibria for the microscopic system, we seek the equilibrium points of the macroscopic system.  From the assumptions $I_0^*\ne 0$ and $I_1^*=0$, the macroscopic equations $\dot{S}=0, \dot{I}_j=0, \dot{R}=0, \dot{D}=0$ ($j=0,1$) become
\begin{eqnarray}
&&\left( \beta_{00} + \beta_{01} \right) I_0^* S=\chi R, \label{nomutanteq:eq:1} \\
&&\left( \beta_{00} S - (\pi_0+\delta_0)\right) I_0^*=0, \nonumber \\  
&&\beta_{01} I_0^* S=0, \nonumber \\
&&\chi R = \pi_0 I_0^*, \label{nomutanteq:eq:4} \\
&&\delta_0 I_0^*=0, \label{nomutanteq:eq:5}
\end{eqnarray}
where it has been taken into account that $I_1^*=0$ implies $\nu^*=1.$ From~\eqref{nomutanteq:eq:5} it is clear that a necessary condition to have an equilibrium of such type is that $\delta_0=0$. In the following, this will be assumed.

To analyze it we discuss two cases: $\chi>0$ and $\chi=0$. Precisely,
\begin{itemize}

\item[(\textit{i})] Case $\chi>0$: In this scenario and using that $I_0^*>0$, equations~\eqref{nomutanteq:eq:1}-\eqref{nomutanteq:eq:4}  become
\begin{equation*}
(\beta_{00} + \beta_{01}) I_0^* S = \chi R, \qquad \beta_{00} S = \pi_0, \qquad \beta_{01} S =0, \qquad \chi R = \pi_0 I_0^*.
\label{nomutant:caseii:chinonneg}
\end{equation*}
From here we have four possible scenarios:
\begin{itemize}
    \item[(\textit{i$_1$})] Case $\beta_{00}=0$ and $\beta_{01}>0$: second and third equation requires $\pi_0$ and $S^*=0$ respectively. Given this result the first and last equation result into $R^*=0$. With this, the equilibrium reads as
    \[
        (S^*,I_0^*,I_1^*,R^*,D^*) = (0, I_0^*, 0,0, 1-I_0^*).
    \]
    \item[(\textit{i$_2$})] Case $\beta_{00}>0$ and $\beta_{01}=0$: second equation leads to $S^*=\pi_0/\beta_{00}$ and last equation to $R^*=\pi_0I_0^*/\chi$. Therefore, the equilibrium under this assumptions is
    \[
        (S^*,I_0^*,I_1^*,R^*,D^*) = \left(\frac{\pi_0}{\beta_{00}}, I_0^*, 0,\frac{\pi_0}{\chi}I_0^*, 1-\frac{\pi_0}{\beta_{00}}-I_0^*\left(1+\frac{\pi_0}{\chi}\right)\right),
    \]
    with $0<I_0^*\leq (1-\pi_0/\beta_{00})/(1+\pi_0/\chi)$ to keep all species between 0 to 1.
    Notice that the condition $\beta_{00}>0$ implies an incompatibility with the microscopic state $\QSmicro^{(1)}$, as it requires $v_0^*=0$.
    \item[(\textit{i$_3$})] Case $\beta_{00}>0$ and $\beta_{01}>0$: to accomplish third equation $S^*=0$ is required. Substituting this into first and second equation is trivial to see that $R=0$ and $\pi_0$ are needed. With this, the equilibrium reads as:
    \[
        (S^*,I_0^*,I_1^*,R^*,D^*) = (0, I_0^*, 0,0, 1-I_0^*).
    \]
    Again, the condition $\beta_{00}>0$ excludes the possibility of having $\QSmicro^{(1)}$ as the micro state at equilibrium.
    \item[(\textit{i$_4$})] Case $\beta_{00}=0$ and $\beta_{01}=0$: second equation implies $\pi_0=0$, fourth equation $R^*=0$ and the other two equation are satisfied no matter the values of $S$ and $I_0^*$. Thus, the equilibrium under this condition is
    \[
        (S,I_0^*,I_1^*,R^*,D) = (S, I_0^*, 0,0, D),
    \]
    with $S$, $I_0^*$ and $D$ such that $S+I_0^*+D=1$.
\end{itemize}

\item[(\textit{ii})] Case $\chi=0$: equations~\eqref{nomutanteq:eq:1}-\eqref{nomutanteq:eq:5} simplify to
\[
(\beta_{00} + \beta_{01})I_0^* S =0, \qquad
(\beta_{00} S - \pi_0)  I_0^* = 0, \qquad
\beta_{01}I_0^*S = 0, \qquad
\pi_0 I_0^* =0.
\]
Taking into account that $I_0^*>0$ it follows that $\pi_0=0$ and the latter system reduces to the following three equations:
\begin{equation*}
(\beta_{00} + \beta_{01}) S =0, \qquad \beta_{00} S =0, \qquad \beta_{01} S =0.
\label{nomutant:case:i}
\end{equation*}
Regarding at the first equation above, it suffices $\beta_{00} + \beta_{01}>0$ to imply $S^*=0$. Otherwise, \textcolor{blue}{if} $\beta_{00}=\beta_{01}=0$ then $S$ remains arbitrary. Thus, the solutions of these two possible scenarios are:
\begin{itemize}
   \item[(\textit{ii$_1$})] Case $\beta_{00}=\beta_{01}=0$:
    \[
        (S,I_0^*,I_1^*,R,D) = (S,I_0^*,0,R,D),
    \]
    with $S$, $I_0^*$, $R$ and $D$ such that $S+I_0^*+R+D=1$.
    \item[(\textit{ii$_2$})] Otherwise, that is, if $\beta_{00} + \beta_{01}>0$, we have:
    \[
        (S^*,I_0^*,I_1^*,R,D) = (0, I_0^*, 0,R, D),
    \]
    with $I_0^*$, $R$ and $D$ such that $I_0^*+R+D=1$. Again, the condition $\beta_{00}>0$
    is incompatible with the microscopic state $\QSmicro^{(1)}$.    
\end{itemize}

\end{itemize}
\qed


\subsection{Proof of Proposition~\ref{prop:coex:eqpoints} (CSE)}
From the assumption $I_0^*\ne 0, I_1^*\ne 0$ and~\eqref{eq:sirs:D} it follows that, necessarily, $\delta_0=\delta_1=0$. That is, both master and mutant variants do not induce mortality. Moreover $0<\nu^*=\nu(I_0^*,I_1^*)<1$ so, in particular, the equilibrium point for the microscopic system is given by Lemma~\ref{lemma:micro}, that is:
\[
\begin{array}{lll}
\QSmicromuc: & \ (g_0^*, g_1^* \, ; v_0^*, v_1^*)  =
\left(  1 - \dfrac{\mu}{\mucrit^*} , \dfrac{\mu}{\mucrit^*} \, ;
\dfrac{\xi_0}{\gamma_0}\, g_0^*, \dfrac{\xi_1}{\gamma_1}\, g_1^*
\right) &\qquad \textrm{if $0 < \mu<\mucrit^*$}  \\[2.5ex]
 & \ (g_0^*, g_1^* \, ; v_0^*, v_1^*)  =
\left( 0, 1 \, ; 0, \dfrac{\xi_1}{\gamma_1} \right) &\qquad \textrm{if $\mucrit^* \leq \mu \leq 1$},
\end{array}
\]
where the critical mutation probability associated to $(I_0^*, I_1^*)$ is defined as
\[
\mucrit^* = 1 - \frac{f_1}{f_0} \left( \frac{1}{\nu^*} -1 \right).
\]

We study the possible equilibria for the macroscopic system. Like in the precedent proofs, we distinguish two cases, according to the common waning immunity rate: (i) $\chi>0$; (ii) $\chi=0$. 

\medskip

\noindent (\textit{i}) Case $\chi>0$. \\[1.1ex] 
The macro equilibrium points come from the solutions of
\begin{equation}
\left( \beta_{00} I_0^* + \beta_{01} I_0^* + \beta_{11}I_1^* \right) S= \chi R, \qquad 
\beta_{00}  S = \pi_0, \qquad 
\left( \beta_{01} I_0^* + \beta_{11}I_1^* \right) S = \pi_1 I_1^*, \qquad 
\pi_0 I_0^* + \pi_1 I_1^* = \chi R,
\label{cc:eqpoints:eq}
\end{equation}
where we have taken into account that $I_1^*>0$. Let us now consider two cases:
\begin{itemize}

\item[(\textit{i$_1$})] Case $\beta_{00}=0$. \\[1.1ex]
From the second expression in~\eqref{cc:eqpoints:eq}, it follows that $\pi_0=0$, and from the fourth one:
\begin{equation*}
R^* = \frac{\pi_1}{\chi} I_1^*.    
\end{equation*}
In order to get $S$, notice that either both $\beta_{01}, \beta_{11}>0$ or $\beta_{01}=\beta_{11}=0$. This is clear from the fact that $0<\nu^*<1$ and, hence,
$\beta_{01}= 0 \Leftrightarrow \left( a_1=0 \ \  \textrm{or} \ \  v_1=0 \right)   \Leftrightarrow \beta_{11}=0.
$
\noindent Therefore, two cases arise:
\begin{itemize}
\item[(\textit{a})] Case $\beta_{01}>0,\beta_{11}>0$.\\[1.1ex]
From equations~\eqref{cc:eqpoints:eq} one gets:
\begin{equation*}
S^* = \frac{\pi_1 I_1^*}{\beta_{01}I_0^* + \beta_{11} I_1^*}.    
\end{equation*}
Thus, the complete equilibrium point is of the form
\[
(S^*,I_0^*, I_1^*, R^*, D^*)   \times (g_0^*, g_1^*; v_0^*, v_1^*) =
\left( \frac{\pi_1 I_1^*}{\beta_{01}^*I_0^* + \beta_{11}^* I_1^*} , I_0^*, I_1^*, 
 \frac{\pi_1}{\chi} I_1^*, D^* \right) \times \QSmicromuc,
\]
provided that $S^* + I_0^* + I_1^* + R^* + D^*=1 $.

\item[(\textit{b})] Case $\beta_{01}=\beta_{11}=0$.\\[1.1ex] 
Having in mind the third equation in~\eqref{cc:eqpoints:eq}, we get that $\pi_1=0$ and so $R^*=0$. There is no constraint for $S$. This means that the macroscopic system equilibrium points are given by $(S, I_0^*, I_1^*, 0, D)$ provided that $S+I_0^*+I_1^*+D=1$. Concerning the restrictions that $\beta_{01}=\beta_{11}=0$ imposes on the microscopic equilibrium point 
$\QSmicromuc$, we have that they come from the fact that $\beta_{01}=\beta_{11}=0 \Leftrightarrow a_1=0$ (and then, $\beta_{01}\equiv0, \beta_{11}\equiv 0$) or $v_1^*=0$ (which in our case is only possible if $\xi_1=0$ or in the limit case $\mu=0$, not included in $\QSmicromuc$). That is, the complete equilibrium points are of the form
\[
(S,I_0^*, I_1^*, R^*, D)   \times (g_0^*, g_1^*; v_0^*, v_1^*) = 
(S, I_0^*, I_1^*, 0, D) \times \QSmicromuc,
\]
with arbitrary $S$, $D$ such that $S+I_0^*+I_1^*+D=1$ and provided $\QSmicromuc$ satisfies $\beta_{01}^*=\beta_{11}^*=0$.
\end{itemize}

\medskip

\item[(\textit{i$_2$})] Case $\beta_{00}>0$. \\[1.1ex]
From~\eqref{cc:eqpoints:eq} it turns out that
\[
S^* = \frac{\pi_0}{\beta_{00}}.
\]
Substituting this value into the third equation in~\eqref{cc:eqpoints:eq}, we have
\[
\frac{\beta_{01}\pi_0}{\beta_{00}} I_0^* + \frac{\beta_{11}}{\beta_{00}}\pi_0 I_1^* = \pi_1 I_1^* 
\]
and, therefore,
\begin{equation}
I_1^* = \frac{\beta_{01}\pi_0}{\beta_{00}\pi_1-\beta_{11}\pi_0} \, I_0^*,
\label{cc:i2}
\end{equation}
provided that 
\[
\pi_1 > \frac{\beta_{11}}{\beta_{00}} \pi_0.
\]
Since $I_1^*>0$, expression~\eqref{cc:i2} implies that $\beta_{01}>0$ and, consequently, $\beta_{11}>0$, where we have taken into account the assertion in case (\textit{i$_1$}) above. Finally, substituting~\eqref{cc:i2} into the fourth equation in~\eqref{cc:eqpoints:eq}, we obtain
\[
R^* = \frac{1}{\chi}\left(\pi_0+\pi_1\frac{\beta_{01}\pi_0}{\beta_{00}\pi_1-\beta_{11}\pi_0}\right) \, I_0^*.
\]
Thus, the complete equilibrium point is of the form
\begin{eqnarray*}
&&{(S^*,I_0^*, I_1^*, R^*, D^*) \times (g_0^*, g_1^*; v_0^*, v_1^*) = }\\[1.3ex] 
&&\qquad  \left( \frac{\pi_0}{\beta_{00}^*}, I_0^*, \frac{\beta_{01}^*\pi_0}{\beta_{00}^*\pi_1-\beta_{11}^*\pi_0} \, I_0^*, \frac{1}{\chi}\left(\pi_0+\pi_1\frac{\beta_{01}^*\pi_0}{\beta_{00}^*\pi_1-\beta_{11}^*\pi_0}\right) \, I_0^*, D^*\right) \times \QSmicromuc
\end{eqnarray*}
\end{itemize}

\medskip

\noindent (\textit{ii}) Case $\chi=0$. \\[1.1ex] 
From equation $\pi_0 I_0^* + \pi_1 I_1^* = \chi R$ we get that $\pi_0=\pi_1=0$ and so, the equations for the macroscopic equilibrium points become
\[
(\beta_{00} I_0^* + \beta_{01} I_0^* + \beta_{11} I_1^*) S=0, \qquad
\beta_{00} I_0^* S =0, \qquad 
( \beta_{01} I_0^* + \beta_{11} I_1^* ) S =0.
\]
As before, we consider two cases:
\begin{itemize}
\item[(\textit{ii$_1$})] Case $\beta_{00}=0$.\\[1.1ex]
As before, we distinguish two possible situations:
\begin{itemize}
\item[(\textit{a})] Case $\beta_{01}>$, $\beta_{11}>0$.\\[1.1ex]
From the equations above it follows that $S^*=0$ and arbitrary $R$. Hence, the complete equilibrium point is of the form
\[
(S^*,I_0^*, I_1^*, R, D)   \times (g_0^*, g_1^*; v_0^*, v_1^*) = 
(0, I_0^*, I_1^*, R, D) \times \QSmicromuc,
\]
with arbitrary $R$, $D$ such that $I_0^*+I_1^*+R+D=1$ and provided $\QSmicromuc$ is compatible with $\beta_{00}^*=0$.

\item[(\textit{b})] Case $\beta_{00}=\beta_{11}=0.$\\[1.1ex]
In this scenario, we have that $S$, $R$, and $D$ are arbitrary. So, the complete equilibrium is of the form
\[
(S,I_0^*, I_1^*, R, D)   \times (g_0^*, g_1^*; v_0^*, v_1^*) = 
(S, I_0^*, I_1^*, R, D) \times \QSmicromuc,
\]
with arbitrary $S$, $R$, and $D$ such that $S+I_0^*+I_1^*+R+D=1$ and provided $\QSmicromuc$ is compatible with the conditions $\beta_{01}^*=\beta_{11}^*=0$.
\end{itemize}

\item[(\textit{ii$_2$})] Case $\beta_{00}>0$.\\[1.1ex]
In this case, it follows that $S=0$. So the complete equilibrium point takes the form
\[
(S^*,I_0^*, I_1^*, R, D)   \times (g_0^*, g_1^*; v_0^*, v_1^*) = 
(0, I_0^*, I_1^*, R, D) \times \QSmicromuc,
\]
with arbitrary $R$, and $D$ such that $I_0^*+I_1^*+R+D=1$
\end{itemize}
\qed


\section{On the Next Generation Matrix method}
\label{app:NGMdetails}
\setcounter{equation}{0}
\renewcommand{\theequation}{\thesection\arabic{equation}}
\setcounter{figure}{0}
The application of the Next Generation Matrix (NGM) method, as 
used in Section~\ref{sec:basic:reproduction:number}, requires the verification of some assumptions on the system (see~\cite{DriesscheWatmough02,DriesscheWatmough08}). Following these references, the system is written in the form
\begin{equation*}
\dot{x}_i=\mathcal{F}_i(x,y)-\mathcal{V}_i(x,y), \qquad
\dot{y}_j=g_j(x,y),
\end{equation*}
where $x$ denotes the infected compartments and $y$ the disease-free compartments. The functions $\mathcal{F}$ and $\mathcal{V}$ and the disease-free subsystem must satisfy several mathematical assumptions, namely:
\begin{itemize}
\item[(\textit{i})] $\mathcal{F}_i(0,y)=0$ and $\mathcal{V}_i(0,y)=0$ $\forall y\geq 0$ and $i=1,\ldots,n$. That is, all the new infections arise solely from secondary infections generated by infected individuals.

\item[(\textit{ii})] $\mathcal{F}_i(x,y)\geq 0$ $\forall x,y\geq 0$ and $i=1,\ldots,n$, reflecting the fact that the rate of new infections cannot be negative.

\item[(\textit{iii})] $ \mathcal{V}_i(x,y)\leq 0$ whenever $x_i=0$, $i=1,\ldots,n$. Since $\mathcal{V}_i$ represents the outflow from compartment $i$, these conditions ensure that the flow is inward whenever the compartment is empty. 

\item[(\textit{iv})] $\sum_{i=1}^n \mathcal{V}_i(x,y) \geq 0$ $\forall x,y\geq 0$, implying that the total outflow from the set of infected compartments is non-negative.

\item[(\textit{v})] The disease-free subsystem $\dot{y}=g(0,y)$ has a unique equilibrium, DFE, which is asymptotically stable.

\end{itemize}

In our particular case, we assume $\chi>0$\footnote{The case $\chi=0$ corresponds to a continuum of disease-free equilibria. In that case, the Van den Driessche--Watmough formulation \cite{DriesscheWatmough02,DriesscheWatmough08} is not satisfied in the strict sense. However, the NGM method can still be applied locally around a chosen disease-free state, usually the fully susceptible one.}. 
It is straightforward to check that decomposition~\eqref{NGM:disease:comp:ode} and~\eqref{NGM:disease-free:comp:ode} in Section~\ref{sec:basic:reproduction:number}
satisfy NGM assumptions (\textit{i})--(\textit{iv}). Concerning (\textit{v}), in the absence of infected individuals, the disease-free subsystem is
\[
\dot{S}=\chi R, \qquad \dot{R}=-\chi R, \qquad \dot{D}=0.
\]
The point $(S,R,D)=(1,0,0)$ is therefore a disease-free equilibrium. To prove its asymptotic stability in the biologically relevant invariant set with $D(0)=0$ and $S(0)+R(0)=1$, observe that $\dot{R}=-\chi R$ implies
$R(t)=\rme^{-\chi t}R(0),$
and hence
$\lim_{t\rightarrow+\infty}R(t)=0.$
Moreover, since $\dot{D}=0$, we have $D(t)=0$ for all $t\geq 0$. Finally, using $\dot{S}=\chi R$,
\[
S(t)=S(0)+\chi\int_0^t R(s)\,ds
= S(0)+R(0)\left(1-\rme^{-\chi t}\right),
\]
so that
\[
\lim_{t\rightarrow+\infty}S(t)=S(0)+R(0)=1.
\]
Therefore, every solution in this disease-free invariant set converges to $(1,0,0)$, which proves assumption (\textit{v}).


\end{document}